\documentclass[12pt]{article}
\usepackage{amsmath,amsfonts,amssymb,amsthm}
\input amssym.def
\begin{document}
\def \Z{\Bbb Z}
\def \C{\Bbb C}
\def \R{\Bbb R}
\def \Q{\Bbb Q}
\def \N{\Bbb N}

\def \A{{\mathcal{A}}}
\def \D{{\mathcal{D}}}
\def \E{{\mathcal{E}}}
\def \L{\mathcal{L}}
\def \S{{\mathcal{S}}}
\def \wt{{\rm wt}}
\def \tr{{\rm tr}}
\def \span{{\rm span}}
\def \Res{{\rm Res}}
\def \Der{{\rm Der}}
\def \End{{\rm End}}
\def \Ind {{\rm Ind}}
\def \Irr {{\rm Irr}}
\def \Aut{{\rm Aut}}
\def \GL{{\rm GL}}
\def \Hom{{\rm Hom}}
\def \mod{{\rm mod}}
\def \ann{{\rm Ann}}
\def \ad{{\rm ad}}
\def \rank{{\rm rank}\;}
\def \<{\langle}
\def \>{\rangle}

\def \g{{\frak{g}}}
\def \h{{\hbar}}
\def \k{{\frak{k}}}
\def \sl{{\frak{sl}}}
\def \gl{{\frak{gl}}}

\def \be{\begin{equation}\label}
\def \ee{\end{equation}}
\def \bex{\begin{example}\label}
\def \eex{\end{example}}
\def \bl{\begin{lem}\label}
\def \el{\end{lem}}
\def \bt{\begin{thm}\label}
\def \et{\end{thm}}
\def \bp{\begin{prop}\label}
\def \ep{\end{prop}}
\def \br{\begin{rem}\label}
\def \er{\end{rem}}
\def \bc{\begin{coro}\label}
\def \ec{\end{coro}}
\def \bd{\begin{de}\label}
\def \ed{\end{de}}

\newcommand{\m}{{\bf m}}
\newcommand{\n}{{\bf n}}
\newcommand{\x}{{\bf x}}
\newcommand{\z}{{\bf z}}
\newcommand{\y}{{\bf y}}
\newcommand{\bft}{{\bf t}}

\newcommand{\nno}{\nonumber}
\newcommand{\nord}{\mbox{\scriptsize ${\circ\atop\circ}$}}
\newtheorem{thm}{Theorem}[section]
\newtheorem{prop}[thm]{Proposition}
\newtheorem{coro}[thm]{Corollary}
\newtheorem{conj}[thm]{Conjecture}
\newtheorem{example}[thm]{Example}
\newtheorem{lem}[thm]{Lemma}
\newtheorem{rem}[thm]{Remark}
\newtheorem{de}[thm]{Definition}
\newtheorem{hy}[thm]{Hypothesis}
\makeatletter
\@addtoreset{equation}{section}
\def\theequation{\thesection.\arabic{equation}}
\makeatother
\makeatletter

\begin{center}
{\Large \bf Toroidal vertex algebras and their modules}
\end{center}

\begin{center}
{Haisheng Li$^{a}$\footnote{Partially supported by NSA grant
H98230-11-1-0161 and China NSF grant (No. 11128103)}, Shaobin Tan$^{b}$\footnote{Partially supported by
China NSF grant (No.10931006) and a grant from the PhD Programs
Foundation of Ministry of Education of China (No.20100121110014).}
and Qing Wang$^{b}$\footnote{Partially supported by China NSF grant
(No.11001229) and the Fundamental Research Funds for the Central
  University (No.2010121003).}\\
$\mbox{}^{a}$Department of Mathematical Sciences\\
Rutgers University, Camden, NJ 08102, USA\\
$\mbox{}^{b}$School of Mathematical Sciences, Xiamen University,
Xiamen 361005, China}
\end{center}

\begin{abstract}
We develop a theory of toroidal vertex algebras and their
modules, and we give a conceptual construction of toroidal
vertex algebras and their modules. As an application, we associate
toroidal vertex algebras and their modules to toroidal Lie
algebras.
\end{abstract}

\section{Introduction}
This is a paper in a series to study extended affine Lie
algebras by using vertex algebra language and techniques. In this
paper, we introduce and study a theory of $(r+1)$-toroidal vertex
algebras with $r$ a positive integer, and we give a general
construction. As an example, we associate $(r+1)$-toroidal vertex
algebras and their modules to $(r+1)$-toroidal Lie algebras.

Toroidal Lie algebras, which are central extensions of multi-loop Lie algebras,
form a special family of what were called extended affine Lie algebras (see \cite{aabgp}),
which are a large family of Lie algebras,
generalizing affine Kac-Moody Lie algebras in a certain natural way.
Affine Kac-Moody Lie algebras are classified as untwisted affine algebras
and twisted affine algebras (see \cite{kac1}), where untwisted affine Lie algebras are
the universal central extensions of $1$-loop Lie algebras while twisted affine
algebras are (or can be realized as) the fixed-point subalgebras of
untwisted affine algebras with respect to Dynkin diagram automorphisms. It is
expected that every extended affine Lie algebra can be realized as the fixed-point subalgebra of
a toroidal Lie algebra with respect to a finite abelian group
(cf. \cite{az}, \cite{aby}, \cite{abfp}, \cite{abp}).

It has been well known (\cite{flm}, \cite{fz}, \cite{dl}; cf. \cite {li-local}, \cite{li-twisted}) 
that vertex algebras can be associated to both twisted
and twisted affine Lie algebras. More specifically, vertex algebras and modules are associated to
highest weight modules for untwisted affine Lie algebras, whereas twisted modules
for those associated vertex algebras are associated to highest weight modules
for the twisted affine Lie algebras.
In the literature, a connection of toroidal Lie algebras with vertex algebras
has also been known (see \cite{bbs}), which uses
one-variable generating functions for toroidal Lie algebras.
On the other hand, it is natural and also we need for various purposes to consider
multi-variable generating functions (cf. \cite{iku}, \cite{ikux}). Then an important question is what kind of
vertex algebra-like structures we can possibly get by using such multi-variable generating functions.
This is one of our motivations for introducing
a theory of toroidal vertex algebras.

In this paper, for a given positive integer $r$ we introduce a notion of $(r+1)$-toroidal vertex algebra.
By definition, an $(r+1)$-toroidal vertex algebra is a vector space $V$ equipped with a linear map
\begin{eqnarray*}
Y(\cdot;x_{0},\x):&& V\rightarrow (\End V)[[x_{0}^{\pm 1},x_{1}^{\pm 1},\dots,x_{r}^{\pm 1}]]\\
&&v\mapsto Y(v;x_{0},\x),
\end{eqnarray*}
where $\x=(x_{1},\dots,x_{r})$,
such that for $u,v\in V$,
$$Y(u;x_{0},\x)v\in V[[x_{1}^{\pm 1},\dots,x_{r}^{\pm 1}]]((x_{0})),$$
and
\begin{eqnarray}\label{ejacobi-rva-introduction}
& &z_{0}^{-1}\delta\left(\frac{x_{0}-y_{0}}{z_{0}}\right)
Y(u;x_{0},\z\y)Y(v;y_{0},\y) -
z_{0}^{-1}\delta\left(\frac{y_{0}-x_{0}}{-z_{0}}\right)
Y(v;y_{0},\y)Y(u;x_{0},\z\y)\nonumber\\
 & &\ \ \ \ \ =
y_{0}^{-1}\delta\left(\frac{x_{0}-z_{0}}{y_{0}}\right)Y(Y(u;z_{0},\z)v;y_{0},\y).
\end{eqnarray}
The existence of a vector ${\bf 1}$, called the vacuum vector, is also assumed, satisfying
$$Y({\bf 1};x_{0},\x)v=v\ \mbox{ and } \ Y(v; x_{0},\x){\bf 1}\in
V[[x_{0},x_{1}^{\pm 1},\dots,x_{r}^{\pm 1}]]\ \ \mbox{ for }v\in
V.$$

Note that we here do not have the full creation property for a vertex algebra $V$:
$$Y(v,x){\bf 1}\in V[[x]]\ \ \mbox{ and }\ \left(Y(v,x){\bf 1}\right)|_{x=0}=v \ \ \mbox{ for }v\in V.$$
In this toroidal vertex algebra theory, the vertex operator map $Y(\cdot;x_{0},\x)$ in general may be {\em not} injective.
That is, the state-field correspondence may be not one-to-one.
Furthermore, for a vertex algebra $V$, one has a canonical operator $D$, defined by
$$D(v)=\left(\frac{d}{dx}Y(v,x){\bf 1}\right)|_{x=0}\ \ \mbox{ for }v\in V,$$
and the following important properties hold:
\begin{eqnarray*}
&&[D,Y(v,x)]=\frac{d}{dx}Y(v,x)\\
&&Y(u,x)v=e^{xD}Y(v,-x)u\ \ \ \mbox{ for }u,v\in V.
\end{eqnarray*}
For a general toroidal vertex algebra, we {\em no longer} have such properties.
These represent the major differences between the notion of toroidal vertex algebra
and that of vertex algebra.

In this paper, we also give a general construction of $(r+1)$-toroidal
vertex algebras. Let $W$ be a general vector space. Set
$$\E(W,r)=\Hom (W,W[[x_{1}^{\pm 1},\dots,x_{r}^{\pm 1}]]((x_{0}))).$$
We consider local subsets $U$ of $\E(W,r)$ in the sense that for any $a(x_{0},\x),b(x_{0},\x)\in U$,
there exists a nonnegative integer $k$ such that
$$(x_{0}-z_{0})^{k}[a(x_{0},\x),b(z_{0},\z)]=0.$$
It is proved that every local subset generates an $(r+1)$-toroidal vertex algebra
in a certain canonical way with $W$ as a module.
Roughly speaking, for $a(x_{0},\x),b(x_{0},\x)\in \E(W,r)$,
we define $a(x_{0},\x)_{m_{0},\m}b(x_{0},\x)\in \E(W,r)$
for $(m_{0},\m)\in \Z\times \Z^{r}$ in terms of
generating function
$$Y_{\E}(a(y_{0},\y);z_{0},\z)b(y_{0},\y)
=\sum_{(m_{0},\m)\in \Z\times \Z^{r}}a(y_{0},y)_{m_{0},\m}b(y_{0},\y)z_{0}^{-m_{0}-1}\z^{-\m}$$
symbolically by
$$Y_{\E}(a(y_{0},\y);z_{0},\z)b(y_{0},\y)=\left(a(x_{0},\x)b(y_{0},\y)\right)|_{x_{0}=y_{0}+z_{0},\x=\y\z}$$
(see Section 3 for the precise definition).
Note that for vertex algebras one uses the usual operator product expansion
$$a(x)b(z)\sim \sum_{n\in \Z}(x-z)^{-n-1}C_{n}(z).$$
For toroidal vertex agebras, we use the following operator product expansion
$$a(x_{0},\x)b(z_{0},\z)\sim \sum_{m_{0}\in \Z,\; \m\in \Z^{r}}
(x_{0}-z_{0})^{-m_{0}-1}\left(\frac{\x}{\z}\right)^{-\m}C_{m_{0},\m}(z_{0},\z).$$

Let $\g$ be a Lie algebra equipped with a symmetric invariant
bilinear form $\<\cdot,\cdot\>$. Associated to
$(\g,\<\cdot,\cdot\>)$, one has a general affine Lie algebra
$$\hat{\g}=\g\otimes \C[t,t^{-1}]\oplus \C {\bf k},$$
which is used to construct vertex algebras. Now, consider a central extension
of the $(r+1)$-loop Lie algebra
$$T_{r+1}(\g)=\g\otimes \C[t_{0}^{\pm 1},t_{1}^{\pm 1},\dots,t_{r}^{\pm 1}]\oplus \C {\bf k},$$
which is referred to as the $(r+1)$-toroidal Lie algebra of $\g$.
We here use this Lie algebra to construct $(r+1)$-toroidal vertex algebras.

Set $B=\g\otimes  \C[t_{0},t_{1}^{\pm 1},\dots,t_{r}^{\pm 1}]+\C {\bf k}$.
Let $\ell$ be a complex number.
Define a $B$-module $(\g+\C)_{\ell}$ where $(\g+\C)_{\ell}=\g\oplus \C$ as a vector space,
with ${\bf k}$ acting as scalar $\ell$, with $\g\otimes  \C[t_{0},t_{1}^{\pm 1},\dots,t_{r}^{\pm 1}]$
acting trivially on $\C$, and with
$$(a\otimes \bft^{\m})\cdot b=[a,b],\ \ \ (a\otimes t_{0}\bft^{\m})\cdot b=\<a,b\>\ell,\ \ \
(a\otimes t_{0}^{n}\bft^{\m})\cdot b=0$$
for $a,b\in \g,\; \m\in \Z^{r},\ n\ge 2$.
Then form an induced module
$$V_{T_{r+1}(\g)}(\ell,0)=U(T_{r+1}(\g))\otimes_{U(B)} (\g+\C)_{\ell},$$
which naturally contains $\g$ as a subspace. Set ${\bf 1}=1\otimes 1\in V_{T_{r+1}(\g)}(\ell,0)$.
It is proved that  there exists a canonical structure of an $(r+1)$-toroidal vertex algebra
on $V_{T_{r+1}(\g)}(\ell,0)$ with $Y(a;x_{0},\x)=a(x_{0},\x)$ for $a\in \g$, where
$$a(x_{0},\x)=\sum_{(m_{0},\m)\in \Z\times \Z^{r}}
\left(a\otimes t_{0}^{m_{0}}\bft^{\m}\right)x_{0}^{-m_{0}-1}\x^{-\m}.$$
Note that $V_{T_{r+1}(\g)}(\ell,0)$ viewed as a $T_{r+1}(\g)$-module
is {\em not} generated by ${\bf 1}$, unlike the situation for the vertex algebras
associated to affine Lie algebras.

Note that in the literature, certain higher dimension analogues of vertex algebras have been studied before
(see \cite{bor-G}, \cite{li-g1}, \cite{nik}, \cite{baknik}), but toroidal vertex algebras are
not vertex algebras in any of those senses. In fact, there is no natural connection
between toroidal Lie algebras and any of those higher dimension analogues.

This paper is organized as follows: In Section 2, we define the
notion of $(r+1)$-toroidal vertex algebra and the notion of module
for an $(r+1)$-toroidal vertex algebra. We also present some basic
results. In Section 3, we give a general construction of
$(r+1)$-toroidal vertex algebras and their modules. In Section 4, we
associate $(r+1)$-toroidal vertex algebras and their modules to
$(r+1)$-toroidal Lie algebras.

\section{Toroidal vertex algebras and their modules}
In this section we define the notion of $(r+1)$-toroidal vertex
algebra and the notion of module for an $(r+1)$-toroidal vertex
algebra with $r$ a positive integer. We present some basic
properties similar to those for ordinary vertex algebras.

For this paper, the scalar field is the field $\C$ of complex
numbers, though it works fine with any field of characteristic zero.
Letters $x,y,z , x_{0},y_{0},x_{0}, x_{1},y_{1},z_{1},\dots$ will be
mutually commuting independent formal variables. Let $r$ be a
positive integer which is fixed throughout this section. For
${\m}=(m_{1},\dots,m_{r})\in \Z^{r}$, set
$${\x}^{\m}=x_{1}^{m_{1}}\cdots x_{r}^{m_{r}}.$$
As a convention we write
$${\x}^{-1}=x_{1}^{-1}\cdots x_{r}^{-1},\ \ \ \
{\x}^{\m-1}=x_{1}^{m_{1}-1}\cdots x_{r}^{m_{r}-1}.$$
We also set
$$\Res_{\x}=\Res_{x_{1}}\cdots \Res_{x_{r}}.$$

\bd{dtva} {\em  An {\em $(r+1)$-toroidal vertex algebra} is a vector
space $V$, equipped with a linear map
\begin{eqnarray*}
Y(\cdot;x_{0},\x): & &V\rightarrow \Hom (V,V[[x_{1}^{\pm
1},\dots,x_{r}^{\pm
1}]]((x_{0})))\subset (\End V)[[x_{0}^{\pm 1},x_{1}^{\pm 1},\dots,x_{r}^{\pm 1}]]\\
& &v\mapsto Y(v;x_{0},\x)=\sum_{(m_{0},\m)\in
\Z^{r+1}}v_{m_{0},\m}x_{0}^{-m_{0}-1}\x^{-\m}
\end{eqnarray*}
and a vector ${\bf 1}\in V$, satisfying the conditions that
$$Y({\bf
1};x_{0},\x)v=v\ \mbox{ and } \ Y(v; x_{0},\x){\bf 1}\in
V[[x_{0},x_{1}^{\pm 1},\dots,x_{r}^{\pm 1}]]\ \ \mbox{ for }v\in
V,$$
 and that for $u,v\in V$,
\begin{eqnarray}\label{ejacobi-rva}
& &z_{0}^{-1}\delta\left(\frac{x_{0}-y_{0}}{z_{0}}\right)
Y(u;x_{0},\z\y)Y(v;y_{0},\y) -
z_{0}^{-1}\delta\left(\frac{y_{0}-x_{0}}{-z_{0}}\right)
Y(v;y_{0},\y)Y(u;x_{0},\z\y)\nonumber\\
 & &\ \ \ \ \ =
y_{0}^{-1}\delta\left(\frac{x_{0}-z_{0}}{y_{0}}\right)Y(Y(u;z_{0},\z)v;y_{0},\y).
\end{eqnarray}} \ed

We also define a notion of {\em $(r+1)$-toroidal
vertex algebra without vacuum}, using all the axioms
above that do not involve the vacuum vector ${\bf 1}$. (Note that a
notion of vertex algebra without vacuum was
introduced and studied by Huang and Lepowsky \cite{hl}.)

\br{rcreation-property}
{\em Note that for a vertex algebra $U$ we have the creation property
$$Y(u,x){\bf 1}\in U[[x]]\ \ \mbox{ and }\ \
\left(Y(u,x){\bf 1}\right)|_{x=0}=u\ \ \mbox{ for }u\in U,$$
which particularly implies that the vertex operator map $Y(\cdot,x)$ is injective.
For an $(r+1)$-toroidal vertex algebra $V$, we do not have this full
creation property (mainly the second part) and in general the map $Y$ may not be injective.}
\er

Let $V$ be an $(r+1)$-toroidal vertex algebra. Just as with ordinary
vertex algebras, applying $\Res_{z_{0}}$ and $\Res_{x_{0}}$ to the
Jacobi identity (\ref{ejacobi-rva}), respectively, we obtain the
following {\em commutator formula} and {\em iterate formula:}
\begin{eqnarray}\label{ecommutator-1}
[Y(u;x_{0},\z\y),Y(v;y_{0},\y)]=\Res_{z_{0}}
y_{0}^{-1}\delta\left(\frac{x_{0}-z_{0}}{y_{0}}\right)Y(Y(u;z_{0},\z)v;y_{0},\y),
\end{eqnarray}
\begin{eqnarray}\label{eiterate}
& &Y(Y(u;z_{0},\z)v;y_{0},\y)=
\Res_{x_{0}}z_{0}^{-1}\delta\left(\frac{x_{0}-y_{0}}{z_{0}}\right)
Y(u;x_{0},\z\y)Y(v;y_{0},\y)\nonumber\\
& &\hspace{3.5cm} -\Res_{x_{0}}
z_{0}^{-1}\delta\left(\frac{y_{0}-x_{0}}{-z_{0}}\right)
Y(v;y_{0},\y)Y(u;x_{0},\z\y) .
\end{eqnarray}

For $v\in V,\; \m\in \Z^{r}$, set
\begin{eqnarray}
Y(v;x_{0},\m)=\sum_{m_{0}\in \Z}v_{m_{0},\m}x_{0}^{-m_{0}-1}
=\Res_{\x}\x^{\m-1}Y(v;x_{0},\x)\in \Hom (V,V((x_{0}))).
\end{eqnarray}
{}From (\ref{ecommutator-1}) we have
\begin{eqnarray}\label{eumv}
[Y(u;x_{0},\m),Y(v;y_{0},\y)]=\Res_{z_{0}}
y_{0}^{-1}\delta\left(\frac{x_{0}-z_{0}}{y_{0}}\right)\y^{\m}Y(Y(u;z_{0},\m)v;y_{0},\y).
\end{eqnarray}
Furthermore,
\begin{eqnarray}
[Y(u;x_{0},\m),Y(v;y_{0},\n)]=\Res_{z_{0}}
y_{0}^{-1}\delta\left(\frac{x_{0}-z_{0}}{y_{0}}\right)Y(Y(u;z_{0},\m)v;y_{0},\m+\n).
\end{eqnarray}
On the other hand, {}from (\ref{eiterate}) we have
\begin{eqnarray}\label{eiterate-comp1}
& &Y(Y(u;z_{0},\m)v;y_{0},\y)=
\Res_{x_{0}}z_{0}^{-1}\delta\left(\frac{x_{0}-y_{0}}{z_{0}}\right)
\y^{-\m}Y(u;x_{0},\m)Y(v;y_{0},\y)\nonumber\\
& &\hspace{2cm} -\Res_{x_{0}}
z_{0}^{-1}\delta\left(\frac{y_{0}-x_{0}}{-z_{0}}\right)
\y^{-\m}Y(v;y_{0},\y)Y(u;x_{0},\m).
\end{eqnarray}

Just as with ordinary vertex algebras, the Jacobi identity axiom is
equivalent to weak commutativity and weak associativity.

\bp{pcomm-assoc} Let $V$ be a vector space equipped with a linear
map $$Y(\cdot; x_{0},\x): V\rightarrow \Hom (V,V[[x_{1}^{\pm 1},\dots, x_{r}^{\pm
1}]]((x_{0}))).$$
For $u,v\in V$, the Jacobi identity
(\ref{ejacobi-rva}) holds if and only if there exists a nonnegative
integer $k$ such that
\begin{eqnarray}\label{efirst-commutativity}
(x_{0}-y_{0})^{k}Y(u;x_{0},\x)Y(v;y_{0},\y)
=(x_{0}-y_{0})^{k}Y(v;y_{0},\y)Y(u;x_{0},\x)
\end{eqnarray}
and
\begin{eqnarray}
z_{0}^{k}Y(Y(u;z_{0},\z)v;y_{0},\y)
=\left((x_{0}-y_{0})^{k}Y(u;x_{0},\z\y)Y(v;y_{0},\y)\right)|_{x_{0}=y_{0}+z_{0}}.
\end{eqnarray}
\ep

Note that the commutativity relation (\ref{efirst-commutativity})
implies
$$(x_{0}-y_{0})^{k}Y(u;x_{0},\x\y)Y(v;y_{0},\y)
\in \Hom \left(W,W[[x_{1}^{\pm 1},y_{1}^{\pm 1},\dots,x_{r}^{\pm 1},y_{r}^{\pm 1}]]((x_{0},y_{0}))\right),$$
so that the substitution
$$\left((x_{0}-y_{0})^{k}Y(u;x_{0},\z\y)Y(v;y_{0},\y)\right)|_{x_{0}=y_{0}+z_{0}}$$
exists in $\Hom \left(W,W[[x_{1}^{\pm 1},y_{1}^{\pm 1},\dots,x_{r}^{\pm 1},y_{r}^{\pm 1}]]
((y_{0}))[[z_{0}]]\right)$.

\br{rtwo-version-wassoc} {\em  A familiar version of weak
associativity is that for any $u,v,w\in V$, there exists a
nonnegative integer $l$ such that
\begin{eqnarray}
(z_{0}+y_{0})^{l}Y(Y(u;z_{0},\z)v;y_{0},\y)w
=(z_{0}+y_{0})^{l}Y(u;z_{0}+y_{0},\z\y)Y(v;y_{0},\y)w.
\end{eqnarray}
The equivalence between the two versions of weak associativity
essentially follows from \cite{ltw} (Lemma 2.9).} \er

For $(r+1)$-toroidal vertex
algebras we have the following skew symmetry:

\bp{pskewsymmetry} Let $V$ be an $(r+1)$-toroidal vertex algebra without vacuum.
For $u,v\in V$, we have
\begin{eqnarray}
Y(Y(u;z_{0},\z)v; y_{0},\y)=e^{z_{0}\frac{\partial}{\partial
y_{0}}}\z^{\y\frac{\partial}{\partial_{\y}}}Y(Y(v;-z_{0},\z^{-1})u;
y_{0},\y),
\end{eqnarray}
where
$$\z^{\y\frac{\partial}{\partial_{\y}}}
=z_{1}^{y_{1}\frac{\partial}{\partial_{y_{1}}}}\cdots z_{r}^{y_{1}\frac{\partial}{\partial_{y_{r}}}}.$$
\ep

\begin{proof} It follows from the standard arguments as in \cite{fhl}.
Notice that starting from the left-hand side of the Jacobi identity
for the pair $(u,v)$, by changing
 $(u,v)$ to $(v,u)$, $(x_{0},y_{0})$ to $(y_{0},x_{0})$, $(z_{0},\z)$ to $(-z_{0},\z^{-1})$,
and then changing $\y$ to $\y\z$, we obtain the left-hand side of
the Jacobi identity for the pair $(v,u)$. In view of this, we have
\begin{eqnarray*}
&&y_{0}^{-1}\delta\left(\frac{x_{0}-z_{0}}{y_{0}}\right)Y(Y(u;z_{0},\z)v;
y_{0},\y)\\
&=&x_{0}^{-1}\delta\left(\frac{y_{0}+z_{0}}{x_{0}}\right)Y(Y(v;-z_{0},\z^{-1})u;
x_{0},\y\z).
\end{eqnarray*}
Applying $\Res_{x_{0}}$ to both sides, we obtain
\begin{eqnarray*}
Y(Y(u;z_{0},\z)v; y_{0},\y)&=&Y(Y(v;-z_{0},\z^{-1})u;
y_{0}+z_{0},\y\z)\nonumber\\
&=&e^{z_{0}\frac{\partial}{\partial
y_{0}}}\z^{\y\frac{\partial}{\partial_{\y}}}Y(Y(v;-z_{0},\z^{-1})u;
y_{0},\y),
\end{eqnarray*}
as desired.
\end{proof}

\bd{dextendedrva} {\em An {\em extended $(r+1)$-toroidal vertex
algebra} is an $(r+1)$-toroidal vertex algebra $V$, equipped with
linear operators $\D_{0},\D_{1},\dots,\D_{r}$, satisfying the
condition that
\begin{eqnarray}\label{eD-etva}
&&\D_{i}({\bf 1})=0\ \ \mbox{ for }0\le i\le r,\nonumber\\
&&[\D_{0},Y(v;x_{0},\x)]=Y(\D_{0}(v);x_{0},\x)=\frac{\partial}{\partial
x_{0}}Y(v;x_{0},\x),\nonumber\\
&&[\D_{j},Y(v;x_{0},\x)]=Y(\D_{j}(v);x_{0},\x)=\left(x_{j}\frac{\partial}{\partial
x_{j}}\right)Y(v;x_{0},\x)
\end{eqnarray}
for $v\in V,\ 1\le j\le r$.} \ed

For an $(r+1)$-toroidal vertex
algebra $V$, we define a {\em derivation} of $V$ to be a linear operator
$D$ on $V$ satisfying the condition that
\begin{eqnarray}
D({\bf 1})=0\ \ \mbox{ and }\ \ [D,Y(v;x_{0},\x)]=Y(D(v);x_{0},\x)\
\ \mbox{ for }v\in V.
\end{eqnarray}
Then the linear operators $\D_{i}\ (0\le i\le r)$ for an extended
$(r+1)$-toroidal vertex algebra are derivations.

{}From Proposition \ref{pskewsymmetry} we immediately have:

\bc{cskew-symmetry} Suppose that $V$ is an extended $(r+1)$-toroidal vertex algebra such that
 $\D_{i}$ $(i=1,\dots,r)$ are semisimple on $V$. Then
\begin{eqnarray}
Y(Y(u;z_{0},\z)v; y_{0},\y)
&=&Y\left(e^{z_{0}\D_{0}}\z^{\bf D}Y(v;-z_{0},\z^{-1})u;
y_{0},\y\right)
\end{eqnarray}
for  $u,v\in V$, where
\begin{eqnarray*}
\z^{\bf D}=z_{1}^{\D_{1}}z_{2}^{\D_{2}}\cdots z_{r}^{\D_{r}}.
\end{eqnarray*}
Furthermore, if the vertex operator map $Y$ is injective, we have
\begin{eqnarray}
Y(u;z_{0},\z)v=e^{z_{0}\D_{0}}\z^{\bf D}Y(v;-z_{0},\z^{-1})u.
\end{eqnarray}
\ec

\bl{lsemisimple}
Let $V$ be an extended $(r+1)$-toroidal vertex algebra with derivations $\D_{0},\D_{1},\dots,\D_{r}$.
Suppose that $v\in V$ is an eigenvector of $\D_{i}$ for $1\le i\le r$ with
eigenvalues $\lambda_{1},\dots,\lambda_{r}$, respectively.
Then $\lambda_{i}\in \Z$ for $i=1,\dots,r$ and
$Y(v;x_{0},\x)\in \x^{\bf \lambda}(\End V)[[x_{0},x_{0}^{-1}]]$.
\el

\begin{proof} From assumption we have
$$\lambda_{i}Y(v;x_{0},\x)=Y(\D_{i}(v);x_{0},\x)=\left(x_{i}\frac{\partial}{\partial
x_{i}}\right)Y(v;x_{0},\x)$$
 for $1\le i\le r$. Notice that
 $$Y(v;x_{0},\x)\in (\End V)[[x_{0}^{\pm 1},x_{1}^{\pm 1},\dots,x_{r}^{\pm 1}]]$$
and that $\left(x_{i}\frac{\partial}{\partial
x_{i}}\right)\x^{\m}=m_{i}\x^{\m}$ for $\m=(m_{1},\dots,m_{r})\in \Z^{r}$.
Then we have $\lambda_{i}\in \Z$ and
$$Y(v;x_{0},\x)=x_{1}^{\lambda_{1}}\cdots x_{r}^{\lambda_{r}}(\End V)[[x_{0},x_{0}^{-1}]],$$
as desired.
\end{proof}

We define notions of toroidal vertex subalgebra, ideal, and
homomorphism, in the obvious ways. The following are straightforward
to prove:

\bl{lsubalgebra-generating} Let $V$ be an $(r+1)$-toroidal vertex
algebra and $U$ a subset of $V$. Denote by $\<U\>$ the linear span
of vectors
$$u^{(1)}_{\m_{1}}\cdots u^{(k)}_{\m_{k}}u$$
for $k\ge 0,\; u^{(1)},\dots,u^{(k)},u\in U\cup \{{\bf 1}\}$,
$\m_{1},\dots,\m_{k}\in \Z^{r+1}$. Then $\<U\>$ is a toroidal vertex
subalgebra of $V$.
\el

\bl{lhomomorphism-generating} Let $V$ and $K$ be $(r+1)$-toroidal
vertex algebras. Suppose that $\psi$ is a linear map from $V$ to
$K$, satisfying $\psi({\bf 1})={\bf 1}$,
\begin{eqnarray}
\psi( Y(u;x_{0},\x)v)=Y(\psi(u);x_{0},\x)\psi(v)
\end{eqnarray}
for $u\in U,\; v\in V$, where $U$ is a generating subset of $V$.
Then $\psi$ is a homomorphism. \el

\bl{lextended-generating} Let $V$ be an $(r+1)$-toroidal vertex
algebra equipped with linear operators $\D_{0},\D_{1},\dots,\D_{r}$,
such that (\ref{eD-etva}) holds for $v\in U$, where $U$ is a
generating subset of $V$. Then $V$ is an extended $(r+1)$-toroidal
vertex algebra.
 \el

For any vertex algebra $V$, the left ideal of generated by the vacuum vector
is the whole vertex algebra $V$, but for high dimension analogs,
this is not true in general. We next study the left ideal generated by ${\bf 1}$ of an $(r+1)$-toroidal vertex algebra.

\bd{dleft-ideal}
{\em Let $V$ be an $(r+1)$-toroidal vertex algebra.
Define $V^{0}$ to be the left ideal generated by ${\bf 1}$,
that is, $V^{0}$ is the subspace linearly spanned by vectors
$$u^{(1)}_{m^{(1)}_{0},\m^{(1)}}\cdots
u^{(k)}_{m^{(k)}_{0},\m^{(k)}}{\bf 1}$$ for $k\ge 0,\; u^{(i)}\in
V,\; (m^{(i)}_{0},\m^{(i)})\in \Z\times \Z^{r}$. }
\ed

First, we have:

\bl{lbasic-fact} Let $V$ be an $(r+1)$-toroidal vertex algebra. For
$u\in V,\ k\in \N,\ \m\in \Z^{r}$, we have
\begin{eqnarray}\label{eright-identity}
&&Y(u_{k,\m}{\bf 1};x_{0},\x)=0,\nonumber\\
&&Y(u_{-k-1,\m}{\bf 1};x_{0},\x)
=\frac{1}{k!}\left(\frac{\partial}{\partial
x_{0}}\right)^{k}Y(u;x_{0},\m)\x^{-\m}.
\end{eqnarray}
 Furthermore,
\begin{eqnarray}\label{eneed}
Y(u_{m_{0},\m}{\bf 1};x_{0},\x)\in \x^{-\m}(\Hom ( V,V((x_{0}))))
\end{eqnarray}
for $m_{0}\in \Z$ and
\begin{eqnarray}\label{esum-sequence}
Y(u;x_{0},\x)=\sum_{\m\in \Z^{r}}Y(u_{-1,\m}{\bf 1};x_{0},\x).
\end{eqnarray}
\el

\begin{proof} Taking $v={\bf 1}$ in (\ref{eiterate}) we get
\begin{eqnarray}\label{eiterate--vacuum}
Y(Y(u;z_{0},\z){\bf 1};y_{0},\y)&=&
\Res_{x_{0}}x_{0}^{-1}\delta\left(\frac{y_{0}+z_{0}}{x_{0}}\right)
Y(u;x_{0},\y\z)\nonumber\\
&=&Y(u;y_{0}+z_{0},\y\z)\nonumber\\
&=&e^{z_{0}\frac{\partial}{\partial y_{0}}}Y(u;y_{0},\y\z).
\end{eqnarray}
{}From this we obtain
\begin{eqnarray}\label{eiterate-comp-vacuum}
Y(Y(u;z_{0},\m){\bf 1};y_{0},\y)
=\y^{-\m}e^{z_{0}\frac{\partial}{\partial y_{0}}}Y(u;y_{0},\m),
\end{eqnarray}
which implies (\ref{eright-identity}).
By (\ref{eiterate-comp-vacuum}), we have
\begin{eqnarray}\label{e-mid-step}
Y(u_{-1,\m}{\bf 1};x_{0},\x) =\x^{-\m}Y(u;x_{0},\m).
\end{eqnarray}
Summing up over $\m\in \Z^{r}$, we obtain (\ref{esum-sequence}).
\end{proof}

The following is the main result about $V^{0}$:

\bp{pva-subalgebra} Let $V$ be an $(r+1)$-toroidal vertex algebra.
 Then $$V^{0}={\rm
span}\{v_{m_{0},\m}{\bf 1}\;|\; v\in V,\; (m_{0}, \m)\in \Z\times
\Z^{r}\},$$ $V^{0}$ is a toroidal vertex subalgebra of $V$, and for
$u\in V^{0}$,
\begin{eqnarray}\label{efinite-property}
Y(u;x_{0},\x)\in \left(\Hom (V, V((x_{0}))\right)[x_{1}^{\pm 1},\dots,x_{r}^{\pm 1}].
\end{eqnarray}
Define a linear map $Y^{0}(\cdot,x_{0}):
V^{0}\rightarrow (\End V)[[x_{0},x_{0}^{-1}]]$ by
$$Y^{0}(u,x_{0})=Y(u;x_{0},\x)|_{\x=1}
\ \ \ \ \mbox{ for }u\in V^{0}.$$ Then $(V^{0},Y^{0},{\bf 1})$
carries the structure of a vertex algebra and $(V,Y^{0})$ carries
the structure of a $V^{0}$-module. Furthermore,  we have
\begin{eqnarray}
Y(u;x_{0},\x)=\sum_{\m\in \Z^{r}}Y^{0}(u_{-1,\m}{\bf
1},x_{0})\x^{-\m}\ \ \ \ \mbox{ for }u\in V.
\end{eqnarray}
\ep

\begin{proof} Set
$$U={\rm span}\{v_{k,\m}{\bf 1}\;|\; v\in V,\; k\in \Z,\; \m\in
\Z^{r}\}\subset V.$$ By definition we have $U\subset V^{0}$. For
$u,v\in V$, applying (\ref{eiterate}) to ${\bf 1}$
and then using the fact that $Y(u;x_{0},\x){\bf 1}\in
V[[x_{0},x_{1}^{\pm 1},\dots,x_{r}^{\pm 1}]]$ we get
\begin{eqnarray}\label{eweak-assoc-vacuum}
Y(u;z_{0}+y_{0},\z\y)Y(v;y_{0},\y){\bf 1}
=Y(Y(u;z_{0},\z)v;y_{0},\y){\bf 1}.
\end{eqnarray}
It follows that $U$ is a left ideal of $V$ and $V^{0}\subset U$.
Therefore $V^{0}=U$. It then follows from (\ref{eneed}) that for any
$v\in V^{0}$,
$$Y(v;x_{0},\x)\in \left(\Hom (V, V((x_{0}))\right)[x_{1}^{\pm 1},\dots,x_{r}^{\pm 1}].$$
{}From this, we see that for every $v\in V^{0}$, $Y^{0}(v,x_{0})$ is
a well defined element in $\Hom (V,V((x_{0})))$. As $V^{0}$ $(=U)$
is a left ideal of $V$, $V^{0}$ is a toroidal vertex subalgebra and
we have
$$Y^{0}(v,x_{0})V^{0}\subset V^{0}((x_{0}))\ \ \mbox{ for }v\in V^{0}.$$

Let $u,v\in V^{0}$. It can be readily seen that the Jacobi identity
(\ref{ejacobi-rva}) for the ordered pair $(u,v)$, after evaluated at
$\y=1,\; \z=1$, gives rise to the Jacobi identity that we need for
$Y^{0}$. Furthermore,  we have  ${\bf 1}={\bf 1}_{-1,{\bf 0}}{\bf
1}\in V^{0}$ and $Y^{0}({\bf 1},x_{0})=Y({\bf
1};x_{0},\x)|_{\x=1}=1$. To show that $(V^{0},Y^{0},{\bf 1})$
carries the structure of a vertex algebra, we need to verify the
creation property. By (\ref{eiterate--vacuum}) we have
$$Y^{0}(Y(u;z_{0},\z){\bf 1},y_{0}){\bf 1}=Y(u;y_{0}+z_{0},\z){\bf
1}.$$ Then $$\lim_{y_{0}\rightarrow 0}Y^{0}(Y(u;z_{0},\z){\bf
1},y_{0}){\bf 1}=Y(u;z_{0},\z){\bf 1}.$$ {}From this it follows that
the creation property holds. Therefore, $(V^{0},Y^{0},{\bf 1})$
carries the structure of a vertex algebra. As it was mentioned
before,
$$Y^{0}(v,x_{0})=Y(v;x_{0},\x)|_{\x=1}\in \Hom (V,V((x_{0})))$$
for $v\in V^{0}$, $Y^{0}({\bf 1},x)=1$, and $Y^{0}$ satisfies the
Jacobi identity for modules for a vertex algebra. Thus $(V,Y^{0})$
carries the structure of a $V^{0}$-module.

Note that for $u\in V,\; (m_{0},\m)\in \Z\times \Z^{r}$, by
(\ref{e-mid-step}) we have
$$Y^{0}(u_{m_{0},\m}{\bf 1},x_{0})=Y(u_{m_{0},\m}{\bf
1};x_{0},\x)|_{\x=1}=Y(u_{m_{0},\m}{\bf
1};x_{0},\m)=\x^{\m}Y(u_{m_{0},\m}{\bf 1};x_{0},\x).$$ Using
(\ref{esum-sequence}) we get
$$Y(u;x_{0},\x)=\sum_{\m\in \Z^{r}}Y(u_{-1,\m}{\bf
1};x_{0},\x)=\sum_{\m\in \Z^{r}}Y^{0}(u_{-1,\m}{\bf
1},x_{0})\x^{-\m},$$ as desired.
\end{proof}

Furthermore, we have:

\bp{pspecial} Let $V$ be an extended $(r+1)$-toroidal
vertex algebra with derivations $\D_{0},\D_{1},\dots,\D_{r}$.
Then $V^{0}$ is an extended $(r+1)$-toroidal vertex subalgebra and
$\D_{j}$ $(1\le j\le r)$ are semisimple on $V^{0}$ with integer eigenvalues.
Furthermore, $\D_{j}$ $(1\le j\le r)$ are derivations of $V^{0}$ viewed as a vertex algebra and
$\D_{0}$ coincides with the $D$-operator $\D$ of the vertex algebra $V^{0}$.
 \ep

\begin{proof} It can be readily seen that the toroidal vertex subalgebra $V^{0}$ is stable under
every derivation of $V$. In particular, $V^{0}$ is stable under
$\D_{0},\D_{1},\dots,\D_{r}$. For $v\in V,\; m_{0}\in \Z,\
\m=(m_{1},\dots,m_{r})\in \Z^{r},\ 1\le j\le r$, we have
\begin{eqnarray}
\D_{j}\left(v_{m_{0},\m}{\bf 1}\right)=v_{m_{0},\m}(\D_{j}{\bf 1})-m_{j}\left(v_{m_{0},\m}{\bf 1}\right)
=-m_{j}\left(v_{m_{0},\m}{\bf 1}\right).
\end{eqnarray}
We see that $\D_{j}$ $(1\le j\le r)$ act semisimply on $V^{0}$ with integer eigenvalues. Furthermore, $\D_{j}$ $(1\le j\le r)$
are derivations of $V^{0}$ viewed as a vertex algebra because
$$[\D_{j},Y^{0}(v,x_{0})]=[\D_{j},Y(v;x_{0},\x)]|_{\x=1}=Y(\D_{j}v;x_{0},\x)|_{\x=1}=Y^{0}(D_{j}v,x_{0})$$
for $v\in V^{0}$.

For $u\in V,\ (m_{0},\m)\in \Z\times \Z^{r}$, we have
\begin{eqnarray*}
&&\lim_{x_{0}\rightarrow 0}\frac{\partial}{\partial x_{0}}Y^{0}(u_{m_{0},\m}{\bf 1},x_{0}){\bf 1}\\
&=&\lim_{x_{0}\rightarrow 0}\left(\frac{\partial}{\partial x_{0}}Y(u_{m_{0},\m}{\bf 1};x_{0},\x){\bf 1}\right)|_{\x=1}\\
&=&\lim_{x_{0}\rightarrow 0}\left(Y(\D_{0}u_{m_{0},\m}{\bf 1};x_{0},\x){\bf 1}\right)|_{\x=1}\\
&=&\lim_{x_{0}\rightarrow 0}Y^{0}(\D_{0}u_{m_{0},\m}{\bf 1},x_{0}){\bf 1}\\
&=&\D_{0}(u_{m_{0},\m}{\bf 1}).
\end{eqnarray*}
Thus $\D_{0}=\D$, the $D$-operator of the vertex algebra $V^{0}$.
\end{proof}

On the other hand, we have:

\bp{precover} Let $U$ be an (ordinary) vertex algebra equipped with
mutually commuting derivations $\D_{1},\dots,\D_{r}$ which act on $U$
semisimply with integer eigenvalues. For $u\in U$, set
$$Y(u;x_{0},\x)=Y(\x^{\bf D}u,x_{0}),$$
where $\x^{\bf D}=x_{1}^{\D_{1}}\cdots x_{r}^{\D_{r}}$.
Then $U$ becomes an extended $(r+1)$-toroidal vertex algebra with $\D_{0}=\D$,
the canonical derivation of vertex algebra $U$. Furthermore, $U^{0}=U$. \ep

\begin{proof} From definition,  for $u\in U$ we have $$Y(u;x_{0},\x)\in (\Hom
(U,U((x_{0}))))[[x_{1}^{\pm 1},\dots,x_{r}^{\pm 1}]].$$
 We also have
$$Y({\bf 1};x_{0},\x)=Y(\x^{\bf D}{\bf 1},x_{0})=Y({\bf 1},x_{0})=1$$
\begin{eqnarray}\label{ecreation-special}
Y(u;x_{0},\x){\bf 1}=Y(\x^{\bf D}u,x_{0}){\bf 1}=e^{x_{0}\D_{0}}\x^{\bf D}u.
\end{eqnarray}
Then the vacuum property holds.  For $u,v\in U$, noticing that
$$Y(Y((\y\z)^{\bf D}u,z_{0})\y^{\bf D}v,y_{0})=Y(\y^{\bf D}Y(\z^{\bf D}u,z_{0})v,y_{0})
=Y(Y(u;z_{0},\z)v;y_{0},\y),$$ we see that the Jacobi identity for a
toroidal vertex algebra holds.

Let $d$ be any derivation of $U$. It is straightforward to show that $d({\bf 1})=0$.
Furthermore, for $u\in U$ we have
$$d\D_{0}(u)=\frac{d}{dx}\left(dY(u,x){\bf 1}\right)|_{x=0}
=\frac{d}{dx}\left(Y(d(u),x){\bf 1}+Y(u,x)d({\bf 1})\right)|_{x=0}=\D_{0}d(u).$$
This proves $d\D_{0}=\D_{0}d$. In particular, we have $\D_{0}\D_{i}=\D_{i}\D_{0}$ for $1\le i\le r$. Using this we get
\begin{eqnarray*}
&&Y(\D_{0}u;x_{0},\x)=Y(\x^{\bf D}\D_{0}u,x_{0})=Y(\D_{0}\x^{\bf D}u,x_{0})
=\frac{\partial}{\partial x_{0}}Y(\x^{\bf D}u,x_{0})\\
&&\ \ \ \ =\frac{\partial}{\partial x_{0}}Y(u;x_{0},\x).
\end{eqnarray*}
For $1\le i\le r$, since $\D_{i}$ acts on $U$
semisimply with integer eigenvalues, we have
$\left(x_{i}\frac{\partial}{\partial x_{i}}\right)\x^{\bf D}u=\D_{i}\x^{\bf D}u$ for $u\in U$. Then
\begin{eqnarray*}
&&Y(\D_{i}u;x_{0},\x)=Y(\x^{\bf D}\D_{i}u,x_{0})=Y(D_{i}\x^{\bf D}u,x_{0})
=x_{i}\frac{\partial}{\partial x_{i}}Y(\x^{\bf D}u,x_{0})\\
&&\ \ \ \ =x_{i}\frac{\partial}{\partial x_{i}}Y(u;x_{0},\x).
\end{eqnarray*}
This proves that $U$ becomes an extended $(r+1)$-toroidal vertex algebra.

Suppose $u\in U$ such that $\D_{j}u=m_{j}u$ for $1\le j\le r$ with $\m=(m_{1},\dots,m_{r})\in \Z^{r}$.
By (\ref{ecreation-special}) we have
 $$Y(u;x_{0},\x){\bf 1}=e^{x_{0}\D_{0}}\x^{\bf D}u=\x^{\m}e^{x_{0}\D_{0}}u,$$
which implies $u=u_{-1,\m}{\bf 1}\in U^{0}$. It follows that $U=U^{0}$.
\end{proof}

Next, we study modules for $(r+1)$-toroidal vertex algebras.

\bd{dmodule} {\em Let $V$ be an $(r+1)$-toroidal vertex algebra. A
{\em $V$-module} is a vector space $W$ equipped with a linear map
\begin{eqnarray*}
Y_{W}(\cdot;x_{0},\x): &&V\rightarrow \Hom (W,W[[x_{1}^{\pm
1},\dots,x_{r}^{\pm
1}]]((x_{0})))\\
&&v\mapsto Y_{W}(v;x_{0},\x),
\end{eqnarray*}
satisfying the condition that $Y_{W}({\bf 1};x_{0},\x)=1_{W}$ and
for $u,v\in V$,
\begin{eqnarray}
&&z_{0}^{-1}\delta\left(\frac{x_{0}-y_{0}}{z_{0}}\right)Y_{W}(u;x_{0},\z\y)Y_{W}(v;y_{0},\y)
\nonumber\\
& &\hspace{2cm}
-z_{0}^{-1}\delta\left(\frac{y_{0}-x_{0}}{-z_{0}}\right)
Y_{W}(v;y_{0},\y)Y_{W}(u;x_{0},\z\y)\nonumber\\
&=&x_{0}^{-1}\delta\left(\frac{y_{0}+z_{0}}{x_{0}}\right)
Y_{W}(Y(u;z_{0},\z)v;y_{0},\y).
\end{eqnarray}}\ed

Just as for an $(r+1)$-toroidal vertex algebra, for a $V$-module, the same commutator
and iterate formulas hold and
the Jacobi identity is equivalent to weak commutativity and weak
associativity. Furthermore, weak associativity follows from the
following {\em iterate formula}
\begin{eqnarray}
& &Y_{W}(Y(u;z_{0},\z)v;y_{0},\y)\nonumber\\
&=&\Res_{x_{0}}z_{0}^{-1}\delta\left(\frac{x_{0}-y_{0}}{z_{0}}\right)
Y_{W}(u;x_{0},\z\y)Y_{W}(v;y_{0},\y)
\nonumber\\
& &\ -\Res_{x_{0}}
z_{0}^{-1}\delta\left(\frac{y_{0}-x_{0}}{-z_{0}}\right)
Y_{W}(v;y_{0},\y)Y_{W}(u;x_{0},\z\y).
\end{eqnarray}
In terms of components, we have
\begin{eqnarray}
& &Y_{W}(u_{m_{0},\m}v;y_{0},\y)=\Res_{x_{0}}\Res_{\x}\x^{\m-1}\y^{-\m}\nonumber\\
& &\cdot
\left((x_{0}-y_{0})^{m_{0}}Y_{W}(u;x_{0},\x)Y_{W}(v;y_{0},\y)
-(-y_{0}+x_{0})^{m_{0}}Y_{W}(v;y_{0},\y)Y_{W}(u;x_{0},\x)\right)\nonumber\\
\end{eqnarray}
for $u,v\in V,\; (m_{0},\m)\in \Z\times \Z^{r}$.

\bd{dmodule-etva} {\em Let $V$ be an extended $(r+1)$-toroidal
vertex algebra. A {\em $V$-module} is a module $(W,Y_{W})$ for $V$
viewed as an $(r+1)$-toroidal vertex algebra, equipped with linear
operators $D_{i}\ (i=0,1,\dots,r)$ on $W$ satisfying that for $v\in
V$,
\begin{eqnarray}\label{emodule-Drelations}
&&[D_{0},Y_{W}(v;x_{0},\x)]=Y_{W}(\D_{0}v;x_{0},\x)=\frac{\partial}{\partial
x_{0}}Y_{W}(v;x_{0},\x),\nonumber\\
&&[D_{i},Y_{W}(v;x_{0},\x)]=Y_{W}(\D_{i}v;x_{0},\x)=x_{i}\frac{\partial}{\partial
x_{i}}Y_{W}(v;x_{0},\x)\ \ \mbox{ for }1\le i\le r.\ \ \ \ \ \
\end{eqnarray}}
\ed

The following is a technical result:

\bl{lextended-generating-module} Let $V$ be an extended
$(r+1)$-toroidal vertex algebra and let $W$ be a module for $V$
viewed as an $(r+1)$-toroidal vertex algebra, equipped with linear
operators $D_{0},D_{1},\dots,D_{r}$, such that
(\ref{emodule-Drelations}) holds for $v\in U$, where $U$ is a
generating subset of $V$. Then $W$ is a $V$-module.
 \el

\begin{proof} We must prove that (\ref{emodule-Drelations})
holds for every $v\in V$. Let $A$ be the collection of $v\in V$ such
that (\ref{emodule-Drelations}) holds. Let $a,b\in A$ and let $w\in
W$. There exists $l\in \N$ such that
$$(x_{0}+z_{0})^{l}Y_{W}(Y(a;z_{0},\z)b;x_{0},\x)w
=(x_{0}+z_{0})^{l}Y_{W}(a;z_{0}+x_{0},\z\x)Y_{W}(b;x_{0},\x)w.$$
Using this we get $a_{m,\n}b\in A$ for $(m,\n)\in \Z\times \Z^{r}$. From
assumption, we have $U\subset A$. Since $U$ generates $V$, we must
have $A=V$. Thus $W$ is a $V$-module.
\end{proof}

\br{rmodule-weakassoc} {\em Note that unlike the case for ordinary
vertex algebras, in the definition of a $V$-module, associativity
alone is not enough because of the lack of the skew symmetry for
toroidal vertex algebras.} \er

\section{Construction of $(r+1)$-toroidal vertex algebras and their modules}
In this section, we give a general construction of $(r+1)$-toroidal
vertex algebras and their modules from local sets of multi-variable
formal vertex operators.

Let $W$ be a vector space and let $r$ be a positive integer, which
are both fixed throughout this section.  Set
$$\E(W,r)=\Hom (W,W[[x_{1}^{\pm 1},\dots,x_{r}^{\pm
1}]]((x_{0})))\subset (\End W)[[x_{0}^{\pm 1},x_{1}^{\pm
1},\dots,x_{r}^{\pm 1}]].$$

Let $a(x_{0},\x),b(x_{0},\x)\in (\End W)[[x_{0}^{\pm 1},x_{1}^{\pm
1},\dots,x_{r}^{\pm 1}]]$.  We say that $a(x_{0},\x)$ and
$b(x_{0},\x)$ are {\em mutually local} if there exists a nonnegative
integer $k$ such that
\begin{eqnarray}\label{elocal-def}
(x_{0}-y_{0})^{k}a(x_{0},\x)b(y_{0},\y)=(x_{0}-y_{0})^{k}b(y_{0},\y)a(x_{0},\x).
\end{eqnarray}
Furthermore, we say that a subset $U$ of $(\End W)[[x_{0}^{\pm
1},x_{1}^{\pm 1},\dots,x_{r}^{\pm 1}]]$ is {\em local} if for any
$a(x_{0},\x),b(x_{0},\x)\in U$, $a(x_{0},\x)$ and $b(x_{0},\x)$ are
mutually local.

Let $a(x_{0},\x),b(x_{0},\x)\in \E(W,r)$. Suppose that $a(x_{0},\x)$
and $b(x_{0},\x)$ are mutually local. Let $k$ be a nonnegative
integer such that (\ref{elocal-def}) holds. Then
\begin{eqnarray}\label{ecompatible}
(x_{0}-y_{0})^{k}a(x_{0},\x)b(y_{0},\y) \in \Hom (W,W[[x_{1}^{\pm
1},\dots,x_{r}^{\pm 1}]]((x_{0},y_{0}))).
\end{eqnarray}
For $a(x_{0},\x),b(x_{0},\x)\in \E(W,r)$, we say that the ordered
pair $(a(x_{0},\x),b(x_{0},\x))$ is {\em compatible} if there exists
a nonnegative integer $k$ such that (\ref{ecompatible}) holds.

\bd{dcomponent-action} {\em Let $a(x_{0},\x),b(x_{0},\x)\in
\E(W,r)$. Assume that $(a(x_{0},\x),b(x_{0},\x))$ is compatible.
Define
$$a(x_{0},\x)_{m_{0},\m}b(x_{0},\x)\in \E(W,r)
 \ \ \ \mbox{ for }(m_{0},\m)\in \Z\times \Z^{r}$$
in terms of generating function
\begin{eqnarray}
Y_{\E}(a(y_{0},\y); z_{0},\z)b(y_{0},\y)=\sum_{(m_{0},\m)\in
\Z\times \Z^{r}} a(y_{0},\y)_{m_{0},\m}b(y_{0},\y)
z_{0}^{-m_{0}-1}\z^{-\m}
\end{eqnarray}
by
\begin{eqnarray}\label{edefinition-yab}
Y_{\E}(a(y_{0},\y);z_{0},\z)b(y_{0},\y) =z_{0}^{-k}
\left((x_{0}-y_{0})^{k}a(x_{0},\z\y)b(y_{0},\y)\right)|_{x_{0}=y_{0}+z_{0}},
\end{eqnarray}
where $k$ is any nonnegative integer such that (\ref{ecompatible})
holds.} \ed

To justify the definition, first note that
$$a(x_{0},\z\y)b(y_{0},\y)\ \mbox{ exists in }(\End W)[[x_{0}^{\pm
1},y_{0}^{\pm 1}, y_{1}^{\pm 1},\dots,y_{r}^{\pm 1},z_{1}^{\pm
1},\dots,z_{r}^{\pm 1}]].$$ Furthermore, with (\ref{ecompatible}) we
have
\begin{eqnarray*}
&&z_{0}^{-k}\left((x_{0}-y_{0})^{k}a(x_{0},\z\y)b(y_{0},\y)\right)|_{x_{0}=y_{0}+z_{0}}\\
&\in& z_{0}^{-k}\Hom (W,W[[y_{1}^{\pm 1},\dots,y_{r}^{\pm
1},z_{1}^{\pm
1},\dots,z_{r}^{\pm 1}]]((y_{0})))[[z_{0}]]\\
&=&\left(\Hom (W,W[[y_{1}^{\pm 1},\dots,y_{r}^{\pm
1}]]((y_{0})))\right)[[z_{1}^{\pm 1},\dots,z_{r}^{\pm 1}]]((z_{0})).
\end{eqnarray*}
Second, it is straightforward to show that the expression on the
right-hand side of (\ref{edefinition-yab}) does not depend on the
choice of $k$ (cf. \cite{li-g1}, \cite{li-ga}).

Set
\begin{eqnarray}
D_{0}=\frac{\partial}{\partial x_{0}}\ \mbox{ and } \
D_{i}=x_{i}\frac{\partial}{\partial x_{i}} \ \ \mbox{  for }1\le
i\le r,
\end{eqnarray}
which act on $\E(W,r)$ in the obvious way. We have:

\bl{ldifferential-operators} Suppose that $a(x_{0},\x),b(x_{0},\x)$
are compatible (resp. local). Then $D_{i}a(x_{0},\x)$ and
$b(x_{0},\x)$ are also compatible (resp. local). Furthermore,
\begin{eqnarray}
[D_{i},Y_{\E}\left(a(x_{0},\x);y_{0},\y\right)]b(x_{0},\x)
=Y_{\E}\left(D_{i}a(x_{0},\x);y_{0},\y\right)b(x_{0},\x)
\end{eqnarray}
for $0\le i\le r$, and
\begin{eqnarray}
&&Y_{\E}\left(D_{0}a(x_{0},\x);y_{0},\y\right)b(x_{0},\x)=\frac{\partial}{\partial
y_{0}}Y_{\E}\left(a(x_{0},\x);y_{0},\y\right)b(x_{0},\x),\label{epartial=0}\\
&&Y_{\E}\left(D_{j}a(x_{0},\x);y_{0},\y\right)b(x_{0},\x)=\left(y_{j}\frac{\partial}{\partial
y_{j}}\right)Y_{\E}\left(a(x_{0},\x);y_{0},\y\right)b(x_{0},\x)\label{epartial=i}
\end{eqnarray}
for $1\le j\le r$. \el

\begin{proof} Let $k$ be a nonnegative integer such that
(\ref{ecompatible}) holds. Then we also have
 \begin{eqnarray}\label{ecompatible+1}
(x_{0}-y_{0})^{k+1}a(x_{0},\x)b(y_{0},\y) \in \Hom (W,W[[x_{1}^{\pm
1},\dots,x_{r}^{\pm 1}]]((x_{0},y_{0}))).
\end{eqnarray}
For $1\le i\le r$, from (\ref{ecompatible}) we immediately have
 \begin{eqnarray*}
(x_{0}-y_{0})^{k}\left(x_{i}\frac{\partial}{\partial
x_{i}}a(x_{0},\x)\right)b(y_{0},\y) \in \Hom (W,W[[x_{1}^{\pm
1},\dots,x_{r}^{\pm 1}]]((x_{0},y_{0}))).
\end{eqnarray*}
For $i=0$, using (\ref{ecompatible+1}) and (\ref{ecompatible}) we
get
 \begin{eqnarray*}
&&(x_{0}-y_{0})^{k+1}(D_{0}a(x_{0},\x))b(y_{0},\y)\\
&=&\frac{\partial}{\partial
x_{0}}\left((x_{0}-y_{0})^{k+1}a(x_{0},\x)b(y_{0},\y)\right)
-(k+1)(x_{0}-y_{0})^{k}a(x_{0},\x)b(y_{0},\y)\\
&\in& \Hom (W,W[[x_{1}^{\pm 1},\dots,x_{r}^{\pm
1}]]((x_{0},y_{0}))).
\end{eqnarray*}
This proves that $D_{i}a(x_{0},\x)$ and $b(x_{0},\x)$ are
compatible. Slightly modifying the arguments, one can show that if
$a(x_{0},\x),b(x_{0},\x)$ are local, then $D_{i}a(x_{0},\x)$ and
$b(x_{0},\x)$ are also local.

The other assertions on $D_{0}$ follow from the same arguments as in
\cite{li-local}. For $D_{i}$ with $1\le i\le r$, notice that
$$\left(x\frac{d}{dx}(x^{n})\right)|_{x=yz}=\left(y\frac{\partial}{\partial
y}\right)(yz)^{n}=\left(z\frac{\partial}{\partial
z}\right)(yz)^{n}$$ for $n\in \Z$. Then the assertions on $D_{i}$
with $1\le i\le r$ (including (\ref{epartial=i})) follow immediately
from definition.
\end{proof}

When $a(x_{0},\x)$ and $b(x_{0},\x)$ are mutually local,
$Y_{\E}(a(y_{0},\y); z_{0},\z)b(y_{0},\y)$ can be defined more
explicitly.

\bl{ldefinition-local} Let $a(x_{0},\x),b(x_{0},\x)\in \E(W,r)$.
Assume that $a(x_{0},\x)$ and $b(x_{0},\x)$ are mutually local. Then
\begin{eqnarray*}
& &Y_{\E}(a(y_{0},\y); z_{0},\z)b(y_{0},\y)\nonumber\\
&=&\Res_{x_{0}}\left(z_{0}^{-1}\delta\left(\frac{x_{0}-y_{0}}{z_{0}}\right)
a(x_{0},\z\y)b(y_{0},\y)
-z_{0}^{-1}\delta\left(\frac{y_{0}-x_{0}}{-z_{0}}\right)
b(y_{0},\y)a(x_{0},\z\y)\right).\
\ \ \
\end{eqnarray*}
In terms of components, we have
\begin{eqnarray} &
&a(y_{0},\y)_{m_{0},\m}b(y_{0},\y)
=\Res_{x_{0}}\Res_{\x}\x^{\m-1}\y^{-\m}\cdot\nonumber\\
& &\hspace{1.5cm}\cdot
\left((x_{0}-y_{0})^{m_{0}}a(x_{0},\x)b(y_{0},\y)
-(-y_{0}+x_{0})^{m_{0}}b(y_{0},\y)a(x_{0},\x)\right) \ \ \ \ \ \
\end{eqnarray}
for $(m_{0},\m)\in \Z\times \Z^{r}$. \el

\begin{proof} Let $k$ be a nonnegative integer such that
$$(x_{0}-y_{0})^{k}a(x_{0},\x)b(y_{0},\y)=(x_{0}-y_{0})^{k}b(y_{0},\y)a(x_{0},\x).$$
Then
$$(x_{0}-y_{0})^{k}a(x_{0},\x)b(y_{0},\y) \in \Hom
(W,W[[\x^{\pm 1},\y^{\pm 1}]]((x_{0},y_{0}))).$$ Using basic
delta-function properties we obtain
\begin{eqnarray*}
&&z_{0}^{k}\Res_{x_{0}}\left(z_{0}^{-1}\delta\left(\frac{x_{0}-y_{0}}{z_{0}}\right)
a(x_{0},\z\y)b(y_{0},\y)
-z_{0}^{-1}\delta\left(\frac{y_{0}-x_{0}}{-z_{0}}\right)
b(y_{0},\y)a(x_{0},\z\y)\right)\\
&=&\Res_{x_{0}}z_{0}^{-1}\delta\left(\frac{x_{0}-y_{0}}{z_{0}}\right)
\left[(x_{0}-y_{0})^{k}a(x_{0},\z\y)b(y_{0},\y]\right)\\
& &\ \ \ \ \
-\Res_{x_{0}}z_{0}^{-1}\delta\left(\frac{y_{0}-x_{0}}{-z_{0}}\right)
\left[(x_{0}-y_{0})^{k}b(y_{0},\y)a(x_{0},\z\y)\right]\\
&=&\Res_{x_{0}}y_{0}^{-1}\delta\left(\frac{x_{0}-z_{0}}{y_{0}}\right)
\left[(x_{0}-y_{0})^{k}a(x_{0},\z\y)b(y_{0},\y)\right]\\
&=&
\left[(x_{0}-y_{0})^{k}a(x_{0},\z\y)b(y_{0},\y)\right]|_{x_{0}=y_{0}+z_{0}}\\
&=&z_{0}^{k}Y_{\E}(a(y_{0},\y); z_{0},\z)b(y_{0},\y),
\end{eqnarray*}
as desired.
\end{proof}

A local subspace $U$ of $\E(W,r)$ is said to be {\em closed} if
$$a(x_{0},\x)_{m_{0},\m}b(x_{0},\x)\in U$$
for all $a(x_{0},\x),b(x_{0},\x)\in U,\; (m_{0},\m)\in \Z\times
\Z^{r}$. As the first main result of this section we have:

\bt{tclosed-local} Let $V$ be a closed local subspace of $\E(W,r)$,
containing $1_{W}$. Then $(V,Y_{\E},1_{W})$ carries the structure of
an $(r+1)$-toroidal vertex algebra and $W$ is a faithful $V$-module
with $Y_{W}(a(x_{0},\x); z_{0},\z)=a(z_{0},\z)$ for $a(x_{0},\x)\in
V$. Furthermore, if $V$ is also stable under $D_{i}$
$(i=0,1,\dots,r)$, then $V$ is an extended $(r+1)$-toroidal vertex
algebra with $\D_{i}=D_{i}$ $(i=0,1,\dots,r)$. \et

\begin{proof} First, for $a(x_{0},\x),b(x_{0},\x)\in V$, there exists
a nonnegative integer $k$ such that
$$(x_{0}-y_{0})^{k}a(x_{0},\x)b(y_{0},\y)
=(x_{0}-y_{0})^{k}b(y_{0},\y)a(x_{0},\x).$$ Then
$a(y_{0},\y)_{m_{0},\m}b(y_{0},\y)\in V$ for $(m_{0},\m)\in \Z\times
\Z^{r}$ and $a(y_{0},\y)_{m_{0},\m}b(y_{0},\y)=0$ for $m_{0}\ge k$.

Second, we check the vacuum property for $1_{W}$. In the definition,
taking $a(x_{0},\x)=1_{W}$, we have
$$1_{m_{0},\m}b(x_{0},\x)=\Res_{y_{0}}\Res_{\y}\y^{\m-1}\x^{-\m}
\left((y_{0}-x_{0})^{m_{0}}-(-x_{0}+y_{0})^{m_{0}}\right)b(x_{0},\x),$$
which equals $b(x_{0},\x)$ when $m_{0}=-1$ and $\m={\bf 0}$ and
equals $0$ otherwise. That is,
\begin{eqnarray}
Y_{\E}(1_{W};z_{0},\z)b(x_{0},\x)=b(x_{0},\x).
\end{eqnarray}
On the other hand, taking $b(x_{0},\x)=1_{W}$, we have
\begin{eqnarray*}
a(x_{0},\x)_{m_{0},\m}1_{W}=\Res_{y_{0}}\Res_{\y}\y^{\m-1}\x^{-\m}
\left((y_{0}-x_{0})^{m_{0}}-(-x_{0}+y_{0})^{m_{0}}\right)a(y_{0},\y),
\end{eqnarray*}
which equals $0$ whenever $m_{0}\ge 0$. For $m_{0}<0$, we have
\begin{eqnarray}
a(x_{0},\x)_{m_{0},\m}1_{W}=\frac{1}{(-m_{0}-1)!}\left(\frac{\partial}{\partial
x_{0}}\right)^{-m_{0}-1}\x^{-\m}a(x_{0},\m).
\end{eqnarray}
Thus
$$\sum_{\m\in \Z^{r}}a(x_{0},\x)_{-1,\m}1_{W}=a(x_{0},\x).$$

Third, we establish the Jacobi identity. Let
$a(x_{0},\x),b(x_{0},\x),c(x_{0},\x)\in V$.  Using definition we
have
\begin{eqnarray}
&&z_{0}^{-1}\delta\left(\frac{x_{0}-y_{0}}{z_{0}}\right)
Y_{\E}(a(t_{0},\bft);x_{0},\z\y)
Y_{\E}(b(t_{0},\bft);y_{0},\y)c(t_{0},\bft)\nonumber\\
&=&z_{0}^{-1}\delta\left(\frac{x_{0}-y_{0}}{z_{0}}\right)
\Res_{s_{0}} x_{0}^{-1}\delta\left(\frac{s_{0}-t_{0}}{x_{0}}\right)
a(s_{0},\z\y\bft)Y_{\E}(b(t_{0},\bft);y_{0},\y)c(t_{0},\bft)\nonumber\\
&&
-z_{0}^{-1}\delta\left(\frac{x_{0}-y_{0}}{z_{0}}\right)\Res_{s_{0}}
x_{0}^{-1}\delta\left(\frac{t_{0}-s_{0}}{-x_{0}}\right)
Y_{\E}(b(t_{0},\bft);y_{0},\y)c(t_{0},\bft)a(s_{0},\z\y\bft)\nonumber\\
&=&\Res_{s_{0}}\Res_{u_{0}}z_{0}^{-1}\delta\left(\frac{x_{0}-y_{0}}{z_{0}}\right)
T,\nonumber
\end{eqnarray}
where
\begin{eqnarray}
T&=&x_{0}^{-1}\delta\left(\frac{s_{0}-t_{0}}{x_{0}}\right)
y_{0}^{-1}\delta\left(\frac{u_{0}-t_{0}}{y_{0}}\right)
a(s_{0},\z\y\bft)b(u_{0},\y\bft)c(t_{0},\bft)\nonumber\\
&&-x_{0}^{-1}\delta\left(\frac{s_{0}-t_{0}}{x_{0}}\right)
y_{0}^{-1}\delta\left(\frac{t_{0}-u_{0}}{-y_{0}}\right)
a(s_{0},\z\y\bft)c(t_{0},\bft)b(u_{0},\y\bft) \nonumber\\
&&-x_{0}^{-1}\delta\left(\frac{t_{0}-s_{0}}{-x_{0}}\right)
y_{0}^{-1}\delta\left(\frac{u_{0}-t_{0}}{y_{0}}\right)
b(u_{0},\y\bft)c(t_{0},\bft)a(s_{0},\z\y\bft)\nonumber\\
&&+x_{0}^{-1}\delta\left(\frac{t_{0}-s_{0}}{-x_{0}}\right)
y_{0}^{-1}\delta\left(\frac{t_{0}-u_{0}}{-y_{0}}\right)
c(t_{0},\bft)b(u_{0},\y\bft)a(s_{0},\z\y\bft).\nonumber
\end{eqnarray}
Let $k$ be a positive integer such that
\begin{eqnarray*}
& &(x_{0}-y_{0})^{k}a(x_{0},\x)b(y_{0},\y)
=(x_{0}-y_{0})^{k}b(y_{0},\y)a(x_{0},\x),\\
& &(x_{0}-y_{0})^{k}a(x_{0},\x)c(y_{0},\y)
=(x_{0}-y_{0})^{k}c(y_{0},\y)a(x_{0},\x),\\
& &(x_{0}-y_{0})^{k}b(x_{0},\x)c(y_{0},\y)
=(x_{0}-y_{0})^{k}c(y_{0},\y)b(x_{0},\x).
\end{eqnarray*}
Using delta-function substitutions we get
\begin{eqnarray*}
& &x_{0}^{k}y_{0}^{k}
z_{0}^{k}z_{0}^{-1}\delta\left(\frac{x_{0}-y_{0}}{z_{0}}\right) T\\
&=&z_{0}^{-1}\delta\left(\frac{x_{0}-y_{0}}{z_{0}}\right)
(s_{0}-t_{0})^{k}(u_{0}-t_{0})^{k}(x_{0}-y_{0})^{k}T\\
&=&z_{0}^{-1}\delta\left(\frac{x_{0}-y_{0}}{z_{0}}\right)
(s_{0}-t_{0})^{k}(u_{0}-t_{0})^{k}(s_{0}-u_{0})^{k} T\\
&=&z_{0}^{-1}\delta\left(\frac{x_{0}-y_{0}}{z_{0}}\right)
t_{0}^{-1}\delta\left(\frac{s_{0}-x_{0}}{t_{0}}\right)
t_{0}^{-1}\delta\left(\frac{u_{0}-y_{0}}{t_{0}}\right)\\
& &\ \ \ \ \ \cdot
\left[(s_{0}-t_{0})^{k}(u_{0}-t_{0})^{k}(s_{0}-u_{0})^{k}
a(s_{0},\z\y\bft)b(u_{0},\y\bft)c(t_{0},\bft)\right].
\end{eqnarray*}

Similarly, we have
\begin{eqnarray}
&&z_{0}^{-1}\delta\left(\frac{y_{0}-x_{0}}{-z_{0}}\right)
Y_{\E}(b(t_{0},\bft);y_{0},\y) Y_{\E}(a(t_{0},\bft);x_{0},\z\y)
c(t_{0},\bft)\nonumber\\
&=&\Res_{s_{0}}\Res_{u_{0}}
z_{0}^{-1}\delta\left(\frac{y_{0}-x_{0}}{-z_{0}}\right)T',\nonumber
\end{eqnarray}
where
\begin{eqnarray}
T'&=&x_{0}^{-1}\delta\left(\frac{s_{0}-t_{0}}{x_{0}}\right)
y_{0}^{-1}\delta\left(\frac{u_{0}-t_{0}}{y_{0}}\right)
b(u_{0},\y\bft)a(s_{0},{\z\y\bft})c(t_{0},{\bft})\nonumber\\
&&-x_{0}^{-1}\delta\left(\frac{t_{0}-s_{0}}{-x_{0}}\right)
y_{0}^{-1}\delta\left(\frac{u_{0}-t_{0}}{y_{0}}\right)
b(u_{0},{\y\bft})c(t_{0},{\bft})a(s_{0},{\z\y\bft})\nonumber\\
&&-x_{0}^{-1}\delta\left(\frac{s_{0}-t_{0}}{x_{0}}\right)
y_{0}^{-1}\delta\left(\frac{t_{0}-u_{0}}{-y_{0}}\right)
a(s_{0},{\z\y\bft})c(t_{0},{\bft})b(u_{0},{\y\bft}) \nonumber\\
&&+x_{0}^{-1}\delta\left(\frac{t_{0}-s_{0}}{-x_{0}}\right)
y_{0}^{-1}\delta\left(\frac{t_{0}-u_{0}}{-y_{0}}\right)
c(t_{0},{\bft})a(s_{0},{\z\y\bft})b(u_{0},{\y\bft}),\nonumber
\end{eqnarray}
and furthermore we have
\begin{eqnarray*}
&&x_{0}^{k}y_{0}^{k}z_{0}^{k}z_{0}^{-1}\delta\left(\frac{y_{0}-x_{0}}{-z_{0}}\right)
Y_{\E}(b(t_{0},\bft);y_{0},\y) Y_{\E}(a(t_{0},\bft);x_{0},\z\y)
c(t_{0},\bft)\nonumber\\
&=&\Res_{s_{0}}\Res_{u_{0}}z_{0}^{-1}\delta\left(\frac{y_{0}-x_{0}}{-z_{0}}\right)
t_{0}^{-1}\delta\left(\frac{s_{0}-x_{0}}{t_{0}}\right)
t_{0}^{-1}\delta\left(\frac{u_{0}-y_{0}}{t_{0}}\right)\\
& &\ \ \ \ \ \cdot
\left[(s_{0}-t_{0})^{k}(u_{0}-t_{0})^{k}(s_{0}-u_{0})^{k}
a(s_{0},\z\y\bft)b(u_{0},\y\bft)c(t_{0},\bft)\right].
\end{eqnarray*}
 On the other hand, we have
\begin{eqnarray}
&&y_{0}^{-1}\delta\left(\frac{x_{0}-z_{0}}{y_{0}}\right)
Y_{\E}(Y_{\E}(a(t_{0},{\bft});z_{0},{\z})
b(t_{0},{\bft});y_{0},{\y})c(t_{0},{\bft})\nonumber\\
&=&y_{0}^{-1}\delta\left(\frac{x_{0}-z_{0}}{y_{0}}\right)
[\Res_{u_{0}}y_{0}^{-1}\delta\left(\frac{u_{0}-t_{0}}{y_{0}}\right)
Y_{\E}(a(u_{0},{\y\bft});z_{0},{\z}) b(u_{0},{\y\bft})c(t_{0},{\bft})\nonumber\\
&&\ \ \ \  -y_{0}^{-1}\delta\left(\frac{t_{0}-u_{0}}{-y_{0}}\right)
c(t_{0},{\bft}) Y_{\E}(a(u_{0},{\y\bft});z_{0},{\z})
b(u_{0},{\y\bft})],\nonumber\\
&=&y_{0}^{-1}\delta\left(\frac{x_{0}-z_{0}}{y_{0}}\right)
\Res_{u_{0}}\Res_{s_{0}}T'',\nonumber
\end{eqnarray}
where
\begin{eqnarray*}
T''&=&y_{0}^{-1}\delta\left(\frac{u_{0}-t_{0}}{y_{0}}\right)
z_{0}^{-1}\delta\left(\frac{s_{0}-u_{0}}{z_{0}}\right)
a(s_{0},{\z\y\bft})b(u_{0},{\y\bft})c(t_{0},{\bft})\nonumber\\
&&-y_{0}^{-1}\delta\left(\frac{u_{0}-t_{0}}{y_{0}}\right)
z_{0}^{-1}\delta\left(\frac{u_{0}-s_{0}}{-z_{0}}\right)
b(u_{0},{\y\bft})a(s_{0},{\z\y\bft})c(t_{0},{\bft})\nonumber\\
&&-y_{0}^{-1}\delta\left(\frac{t_{0}-u_{0}}{-y_{0}}\right)
z_{0}^{-1}\delta\left(\frac{s_{0}-u_{0}}{z_{0}}\right)
c(t_{0},{\bft})a(s_{0},{\z\y\bft})b(u_{0},{\y\bft})\nonumber\\
&&+y_{0}^{-1}\delta\left(\frac{t_{0}-u_{0}}{-y_{0}}\right)
z_{0}^{-1}\delta\left(\frac{u_{0}-s_{0}}{-z_{0}}\right)
c(t_{0},{\bft})b(u_{0},{\y\bft})a(s_{0},{\z\y\bft}),
\end{eqnarray*}
and furthermore we have
\begin{eqnarray*}
&&x_{0}^{k}y_{0}^{k}z_{0}^{k}y_{0}^{-1}\delta\left(\frac{x_{0}-z_{0}}{y_{0}}\right)
Y_{\E}(Y_{\E}(a(t_{0},{\bft});z_{0},{\z})
b(t_{0},{\bft});y_{0},{\y})c(t_{0},{\bft})\nonumber\\
&=&\Res_{s_{0}}\Res_{u_{0}}y_{0}^{-1}\delta\left(\frac{x_{0}-z_{0}}{y_{0}}\right)
t_{0}^{-1}\delta\left(\frac{s_{0}-x_{0}}{t_{0}}\right)
t_{0}^{-1}\delta\left(\frac{u_{0}-y_{0}}{t_{0}}\right)\\
& &\ \ \ \ \ \cdot
\left[(s_{0}-t_{0})^{k}(u_{0}-t_{0})^{k}(s_{0}-u_{0})^{k}
a(s_{0},\z\y\bft)b(u_{0},\y\bft)c(t_{0},\bft)\right].
\end{eqnarray*}
 Note that
\begin{eqnarray*}
&&(x_{0}-y_{0})^{k}x_{0}^{-1}\delta\left(\frac{s_{0}-t_{0}}{x_{0}}\right)
y_{0}^{-1}\delta\left(\frac{u_{0}-t_{0}}{y_{0}}\right)\\
&=&(s_{0}-u_{0})^{k}x_{0}^{-1}\delta\left(\frac{s_{0}-t_{0}}{x_{0}}\right)
y_{0}^{-1}\delta\left(\frac{u_{0}-t_{0}}{y_{0}}\right).
\end{eqnarray*}
Then
\begin{eqnarray*}
&&x_{0}^{k}y_{0}^{k}z_{0}^{k}z_{0}^{-1}\delta\left(\frac{x_{0}-y_{0}}{z_{0}}\right)
Y_{\E}(a(t_{0},\bft);x_{0},\z\y)
Y_{\E}(b(t_{0},\bft);y_{0},\y)c(t_{0},\bft)\nonumber\\
& & -z_{0}^{-1}\delta\left(\frac{y_{0}-x_{0}}{-z_{0}}\right)
Y_{\E}(b(t_{0},\bft);y_{0},\y) Y_{\E}(a(t_{0},\bft);x_{0},\z\y)
c(t_{0},\bft)\nonumber\\
&=&x_{0}^{k}y_{0}^{k}z_{0}^{k}y_{0}^{-1}\delta\left(\frac{x_{0}-z_{0}}{y_{0}}\right)
Y_{\E}(Y_{\E}(a(t_{0},{\bft});z_{0},{\z})
b(t_{0},{\bft});y_{0},{\y})c(t_{0},{\bft}),
\end{eqnarray*}
{}from which we immediately get
\begin{eqnarray*}\label{ejacobi-class}
& &z_{0}^{-1}\delta\left(\frac{x_{0}-y_{0}}{z_{0}}\right)
Y_{\E}(a;x_{0},\z\y)Y_{\E}(b;y_{0},\y)c
\nonumber\\
& &\hspace{2cm} -
z_{0}^{-1}\delta\left(\frac{y_{0}-x_{0}}{-z_{0}}\right)
Y_{\E}(b;y_{0},\y)Y_{\E}(a;x_{0},\z\y)c\nonumber\\
&=&y_{0}^{-1}\delta\left(\frac{x_{0}-z_{0}}{y_{0}}\right)
Y_{\E}(Y_{\E}(a;z_{0},\z)b;y_{0},\y)c,
\end{eqnarray*}
as desired. This establishes Jacobi identity. Therefore,
$(V,Y_{\E},1_{W})$ carries the structure of an $(r+1)$-toroidal
vertex algebra.

For $a(x_{0},\x)\in V$, set $Y_{W}(a(x_{0},\x);
z_{0},\z)=a(z_{0},\z)$. Then, for $a(x_{0},\x),b(x_{0},\x)\in V$,
$w\in W$, we have
\begin{eqnarray*}
&&Y_{W}\left(Y_{\E}(a(x_{0},\x);z_{0},\z)b(x_{0},\x);y_{0},\y\right)w\\
&=&\left(Y_{\E}(a(x_{0},\x);z_{0},\z)b(x_{0},\x)\right)|_{x_{0}=y_{0},\x=\y}w\\
&=&\left(Y_{\E}(a(y_{0},\y);z_{0},\z)b(y_{0},\y)\right)w\\
&=&\Res_{t_{0}}z_{0}^{-1}\delta\left(\frac{t_{0}-y_{0}}{z_{0}}\right)
a(t_{0},\z\y)b(y_{0},\y)w\\
&&-\Res_{t_{0}}z_{0}^{-1}\delta\left(\frac{y_{0}-t_{0}}{-z_{0}}\right)
b(y_{0},\y)a(t_{0},\z\y)w\\
&=&\Res_{t_{0}}z_{0}^{-1}\delta\left(\frac{t_{0}-y_{0}}{z_{0}}\right)
Y_{W}(a(x_{0},\x);t_{0},\z\y)Y_{W}(b(x_{0},\x);y_{0},\y)w\\
&&-\Res_{t_{0}}z_{0}^{-1}\delta\left(\frac{y_{0}-t_{0}}{-z_{0}}\right)
Y_{W}(b(x_{0},\x);y_{0},\y)Y_{W}(a(x_{0},\x);t_{0},\z\y)w.
\end{eqnarray*}
Therefore, $(W,Y_{W})$ carries the structure of a $V$-module which
is faithful as $V\subset \E(W,r)$. For the last assertion, it is
clear from Lemma \ref{ldifferential-operators}.
\end{proof}

Next, we shall show that each local subset of $\E(W,r)$ gives rise
to a closed local subspace. The following technical result follows
from the same arguments as in the single variable case (see
\cite{li-local}):

\bl{linduction} Assume that elements
$a(x_{0},\x),b(x_{0},\x),c(x_{0},\x)$ of $\E(W,r)$ are pair-wise
local. Then for any $(m_{0},\m)\in \Z\times \Z^{r}$,
$a(x_{0},\x)_{m_{0},\m}b(x_{0},\x)$ and $c(x_{0},\x)$ are local. \el

We have:

\bp{pmaximal-tva} Let $V$ be a maximal local subspace of $\E(W,r)$.
Then $V$ contains $1_{W}$, is closed and stable under
$D_{0},D_{1},\dots,D_{r}$. Furthermore, $(V,Y_{\E},1_{W})$ carries
the structure of an extended $(r+1)$-toroidal vertex algebra. \ep

\begin{proof} Notice that $V+\C 1_{W}$ is local and contains $V$.
With $V$ maximal we have $V+\C 1_{W}=V$. Thus $1_{W}\in V$. Let
$a(x_{0},\x),b(x_{0},\x)\in V,\; (m_{0},\m)\in \Z\times \Z^{r}$. By
Lemma \ref{linduction}, for any $c(x_{0},\x)\in V$,
$a(x_{0},\x)_{m_{0},\m}b(x_{0},\x)$ and $c(x_{0},\x)$ are local. In
particular, $a(x_{0},\x)_{m_{0},\m}b(x_{0},\x)$ is local with
$a(x_{0},\x)$ and $b(x_{0},\x)$. Furthermore, using Lemma
\ref{linduction} again, we see that
$a(x_{0},\x)_{m_{0},\m}b(x_{0},\x)$ is local with itself. Thus $V+\C
a(x_{0},\x)_{m_{0},\m}b(x_{0},\x)$ is local. With $V$ maximal we
must have $a(x_{0},\x)_{m_{0},\m}b(x_{0},\x)\in V$. This proves that
$V$ is closed. Using a similar argument and using Lemma
\ref{ldifferential-operators} (the first part) we see that $V$ is
stable under $D_{0},D_{1},\dots,D_{r}$. Then it follows from Theorem
\ref{tclosed-local} that $V$ is an extended $(r+1)$-toroidal vertex
algebra.
\end{proof}

Furthermore, we have:

\bt{tclosedness-tva} Let $U$ be a local subset of $\E(W,r)$. Denote
by $\<U\>$ the linear span of the elements of the form
$$a^{(1)}(x_{0},\x)_{m_{0}^{(1)},\m^{(1)}}\cdots
a^{(k)}(x_{0},\x)_{m_{0}^{(k)},\m^{(k)}}b(x_{0},\x)$$ with $k\in
\N,\; a^{(i)}(x_{0},\x),b(x_{0},\x)\in U\cup \{1_{W}\},\;
(m_{0}^{(i)},\m^{(i)})\in \Z\times \Z^{r}$. Then $\<U\>$ is the
(unique) smallest closed local subspace containing $U\cup
\{1_{W}\}$, and $(\<U\>,Y_{\E},1_{W})$ carries the structure of an
$(r+1)$-toroidal vertex algebra with $W$ as a faithful module.
Furthermore, $\C[D_{0},D_{1},\dots,D_{r}]\<U\>$ is an extended
$(r+1)$-toroidal vertex algebra and $W$ is a module for
$\C[D_{0},D_{1},\dots,D_{r}]\<U\>$ viewed as an $(r+1)$-toroidal
vertex algebra. \et

\begin{proof} By Zorn's Lemma, there exists a maximal local subspace
$V$ of $\E(W,r)$, containing $U$. By Proposition \ref{pmaximal-tva},
$V$ contains $1_{W}$ and is closed, and $V$ is an extended
$(r+1)$-toroidal vertex algebra with $\D_{i}=D_{i}$
$(i=0,1,\dots,r)$. As $U\subset V$, we have $\<U\>\subset V$, so
that $\<U\>$ is local. It follows from induction and
(\ref{eiterate}) that $\<U\>$ is closed. Then $\<U\>$ is an
$(r+1)$-toroidal vertex algebra with $W$ as a faithful module. On
the other hand, we have
$$\C[D_{0},D_{1},\dots,D_{r}]\<U\>\subset V.$$
It follows from Lemma \ref{ldifferential-operators} and induction
that $\C[D_{0},D_{1},\dots,D_{r}]\<U\>$ is a subalgebra of $V$.
Consequently, $\C[D_{0},D_{1},\dots,D_{r}]\<U\>$ is an extended
$(r+1)$-toroidal vertex algebra.
\end{proof}

We shall also need the following result:

\bl{ltransfer} Let $V$ be a closed local subspace of $\E(W,r)$, let
$$a(x_{0},\x),b(x_{0},\x), c_{0}(x_{0},\x),c_{1}(x_{0},\x),\dots,
c_{k}(x_{0},\x)\in V,$$ and let $\m\in \Z^{r}$. If
\begin{eqnarray}
[a(x_{0},\m),b(y_{0},\y)]=\y^{\m}\sum_{j=0}^{k}c_{j}(y_{0},\y)
\frac{1}{j!}\left(\frac{\partial}{\partial y_{0}}\right)^{j}
x_{0}^{-1}\delta\left(\frac{y_{0}}{x_{0}}\right),
\end{eqnarray}
where $a(x_{0},\x)=\sum_{\m\in \Z^{r}}a(x_{0},\m)\x^{-\m}$, then
\begin{eqnarray}
a(x_{0},\x)_{j,\m}b(x_{0},\x)=c_{j}(x_{0},\x)
\end{eqnarray}
for $0\le j\le k$ and $a(x_{0},\x)_{j,\m}b(x_{0},\x)=0$ for $j>k$.
 \el

\begin{proof} For $j\ge 0,\; \m\in \Z^{r}$, from Lemma \ref{ldefinition-local} we have
\begin{eqnarray*}
& &a(y_{0},\y)_{j,\m}b(y_{0},\y)\\
&=&\Res_{x_{0}}\Res_{\x}\x^{\m-1}\y^{-\m}(x_{0}-y_{0})^{j}[a(x_{0},\x),b(y_{0},\y)]\\
&=&\Res_{x_{0}}(x_{0}-y_{0})^{j}\y^{-\m}[a(x_{0},\m),b(y_{0},\y)]\\
&=&\Res_{x_{0}}(x_{0}-y_{0})^{j}\sum_{i=0}^{k}c_{i}(y_{0},\y)
\frac{1}{i!}\left(\frac{\partial}{\partial y_{0}}\right)^{i}
y_{0}^{-1}\delta\left(\frac{x_{0}}{y_{0}}\right)\\
&=&\Res_{x_{0}}(x_{0}-y_{0})^{j}\sum_{i=0}^{k}c_{i}(y_{0},\y) \left(
(x_{0}-y_{0})^{-i-1}-(-y_{0}+x_{0})^{-i-1}\right).
\end{eqnarray*}
If $j>k$, we have $(x_{0}-y_{0})^{j}\left(
(x_{0}-y_{0})^{-i-1}-(-y_{0}+x_{0})^{-i-1}\right)=0$ for $0\le i\le
k$, so that $a(y_{0},\y)_{j,\m}b(y_{0},\y)=0$. If $0\le j\le i$, we
have
\begin{eqnarray*}
&&\Res_{x_{0}}(x_{0}-y_{0})^{j}\left((x_{0}-y_{0})^{-i-1}-(-y_{0}+x_{0})^{-i-1}\right)\\
&=&\Res_{x_{0}}\left((x_{0}-y_{0})^{-(i-j)-1}-(-y_{0}+x_{0})^{-(i-j)-1}\right)\\
&=&\delta_{ij}.
\end{eqnarray*}
In this case, we get
\begin{eqnarray*}
a(y_{0},\y)_{j,\m}b(y_{0},\y)=c_{j}(y_{0},\y),
\end{eqnarray*}
as desired.
\end{proof}

Furthermore we have:

\bp{ptransfer} Suppose that $V$ is a closed local subspace of
$\E(W,r)$ and let
$$a(x_{0},\x),b(x_{0},\x), c_{0}(x_{0},\x),c_{1}(x_{0},\x),\dots,
c_{k}(x_{0},\x)\in V,\;\m\in \Z^{r}.$$ If
\begin{eqnarray}
[a(x_{0},\m),b(y_{0},\y)]=\y^{\m}\sum_{j=0}^{k}c_{j}(y_{0},\y)
\frac{1}{j!}\left(\frac{\partial}{\partial
y_{0}}\right)^{j} x_{0}^{-1}\delta\left(\frac{y_{0}}{x_{0}}\right),
\end{eqnarray}
where $a(x_{0},\x)=\sum_{\m\in \Z^{r}}a(x_{0},\m)\x^{-\m}$, then
\begin{eqnarray}
&&[Y_{\E}(a(t_{0},{\bft});x_{0},\m),Y_{\E}(b(t_{0},\bft);y_{0},\y)]\nonumber\\
&=&\y^{\m}\sum_{j=0}^{k}Y_{\E}(c_{j}(t_{0},\bft);y_{0},\y)
\frac{1}{j!}\left(\frac{\partial}{\partial y_{0}}\right)^{j}
x_{0}^{-1}\delta\left(\frac{y_{0}}{x_{0}}\right).
\end{eqnarray}
In terms of components,
\begin{eqnarray}
&&[Y_{\E}(a(t_{0},{\bft});x_{0},\m),Y_{\E}(b(t_{0},\bft);y_{0},\n)]\nonumber\\
&=&\sum_{j=0}^{k}Y_{\E}(c_{j}(t_{0},\bft);y_{0},\m+\n)
\frac{1}{j!}\left(\frac{\partial}{\partial y_{0}}\right)^{j}
x_{0}^{-1}\delta\left(\frac{y_{0}}{x_{0}}\right).
\end{eqnarray}
\ep

\begin{proof} We may assume that $V$ contains $1_{W}$, so that $V$ is
an $(r+1)$-toroidal vertex algebra. By Lemma \ref{ltransfer} we have
\begin{eqnarray*}
a(x_{0},\x)_{j,\m}b(x_{0},\x)=c_{j}(x_{0},\x)
\end{eqnarray*}
for $0\le j\le k$ and $a(x_{0},\x)_{j,\m}b(x_{0},\x)=0$ for $j>k$.
Then it follows from (\ref{eumv}) immediately.
\end{proof}

As we need, next we establish an analog of a theorem of Xu (see
\cite{xu}; cf. \cite{ll}).

\bt{tanalogue} Let $V$ be a vector space equipped with a linear map
\begin{eqnarray*}
Y(\cdot;x_{0},\x):&& V\rightarrow \Hom (V,V[[x_{1}^{\pm
1},\dots,x_{r}^{\pm 1}]]((x_{0}))),\\
&&v\mapsto Y(v;x_{0},\x)=\sum_{m\in \Z,\n\in \Z^{r}}
v_{m,\n}x_{0}^{-m-1}\x^{-n}
\end{eqnarray*}
and let $U$ be a subset of $V$, satisfying the conditions that
$\{Y(u;x_{0},\x)\ |\ u\in U\}$ is a local subset of $\E(V,r)$, $V$
is linearly spanned by vectors
$$u^{(1)}_{\m_{1}}\cdots u^{(k)}_{\m_{k}}u$$
for $u^{(1)},\dots,u^{(k)},\; u\in U,\;  \m_{1},\dots,\m_{k}\in
\Z^{r+1}$ with $k\in \N$, and that
\begin{eqnarray}\label{eiterate-thm}
& &Y(Y(u;z_{0},\z)v;y_{0},\y)=
\Res_{x_{0}}z_{0}^{-1}\delta\left(\frac{x_{0}-y_{0}}{z_{0}}\right)
Y(u;x_{0},\z\y)Y(v;y_{0},\y)\nonumber\\
& &\hspace{2.5cm} -\Res_{x_{0}}
z_{0}^{-1}\delta\left(\frac{y_{0}-x_{0}}{-z_{0}}\right)
Y(v;y_{0},\y)Y(u;x_{0},\z\y)
\end{eqnarray}
for $u\in U,\; v\in V$. Then $(V,Y)$ carries the structure of an
$(r+1)$-toroidal vertex algebra without vacuum. \et

\begin{proof} Set
$$\tilde{U}=\{ Y(u;x_{0},\x)\;|\; u\in U\}\subset \E(V,r).$$
By assumption, $\tilde{U}$ is a local subset. Then by Theorem
\ref{tclosedness-tva}, $\tilde{U}$ generates an $(r+1)$-toroidal
vertex algebra $\<\tilde{U}\>$, which is also a local subspace. For
$a\in U,\; v\in V$, combining (\ref{eiterate-thm}) with Lemma
\ref{ldefinition-local} we get
\begin{eqnarray}\label{eu-property}
&&Y(Y(a;z_{0},\z)v;y_{0},\y)\nonumber\\
&=&\Res_{x_{0}}z_{0}^{-1}\delta\left(\frac{x_{0}-y_{0}}{z_{0}}\right)
Y(a;x_{0},\z\y)Y(v;y_{0},\y)\nonumber\\
& &\hspace{0.5cm} -\Res_{x_{0}}
z_{0}^{-1}\delta\left(\frac{y_{0}-x_{0}}{-z_{0}}\right)
Y(v;y_{0},\y)Y(a;x_{0},\z\y)\nonumber\\
 &=&Y_{\E}(Y(a;y_{0},\y); z_{0},\z)Y(v;y_{0},\y).
 \end{eqnarray}
Then it follows from the span assumption that $Y$ maps $V$ into
$\<\tilde{U}\>$. Since $\<\tilde{U}\>$ is a local subspace of
$\E(W,r)$, $\{ Y(v;x_{0},\x)\;|\; v\in V\}$ as a subset of
$\<\tilde{U}\>$ is local. Now it remains to establish weak
associativity, which will be achieved in the following by using
iterate formula and induction.

Let $K$ consist of each $u\in V$ such that (\ref{eiterate-thm})
holds for every $v\in V$. Now, let $a,b\in K,\ v\in V$. From
(\ref{eiterate-thm}) with $(a,b)$ in the places of $(u,v)$, there
exists a nonnegative integer $l$ such that
\begin{eqnarray*}
&&(x_{0}+y_{0})^{l}Y(Y(a;x_{0},\x)b;y_{0},\y)v\\
&=&(x_{0}+y_{0})^{l}Y(a;x_{0}+y_{0},\x\y)Y(b;y_{0},\y)v.
\end{eqnarray*}
Using this, (\ref{eu-property}), and the weak associativity for
$\<\tilde{U}\>$,  replacing $l$ with a large one if necessary, we
have
\begin{eqnarray*}
&&(x_{0}+y_{0})^{l}Y(Y(Y(a;x_{0},\x)b;y_{0},\y)v;z_{0},\z)\\
&=&(x_{0}+y_{0})^{l}Y\left(Y(a;x_{0}+y_{0},\x\y)Y(b;y_{0},\y)v;z_{0},\z\right)\\
&=&(x_{0}+y_{0})^{l}Y_{\E}\left(\overline{Y(a;z_{0},\z)};x_{0}+y_{0},\x\y\right)
\overline{Y\left(Y(b;y_{0},\y)v;z_{0},\z\right)}\\
&=&(x_{0}+y_{0})^{l}Y_{\E}\left(\overline{Y(a;z_{0},\z)};x_{0}+y_{0},\x\y\right)
Y_{\E}\left(\overline{Y(b;z_{0},\z)};y_{0},\y\right)\overline{Y(v;z_{0},\z)}\\
&=&(x_{0}+y_{0})^{l}Y_{\E}\left(Y_{\E}\left(\overline{Y(a;z_{0},\z)};x_{0},\x\right)
\overline{Y(b;z_{0},\z)};y_{0},\y\right)\overline{Y(v;z_{0},\z)}\\
&=&(x_{0}+y_{0})^{l}
Y_{\E}\left(\overline{Y\left(Y(a;x_{0},\x)b;z_{0},\z\right)};y_{0},\y\right)
\overline{Y(v;z_{0},\z)},
 \end{eqnarray*}
where $\overline{X}=X$ for all the bar objects; the only purpose is to make the equation easier to read.
Multiplying both sides by formal series $(y_{0}+x_{0})^{-l}$ we get
\begin{eqnarray*}
&&Y(Y(Y(a;x_{0},\x)b;y_{0},\y)v;z_{0},\z)\\
&=&Y_{\E}\left(\overline{Y\left(Y(a;x_{0},\x)b;z_{0},\z\right)};y_{0},\y\right)
\overline{Y(v;z_{0},\z)}.
 \end{eqnarray*}
Furthermore, using  Lemma \ref{ldefinition-local}, we get
\begin{eqnarray*}
&&Y(Y(Y(a,x_{0},\x)b;y_{0},\y)v;z_{0},\z)\\
&=&Y_{\E}\left(Y(Y(a,x_{0},\x)b;z_{0},\z);y_{0},\y\right)
Y(v;z_{0},\z)\\
&=&\Res_{x_{1}}y_{0}^{-1}\delta\left(\frac{x_{1}-z_{0}}{y_{0}}\right)
Y(Y(a;x_{0},\x)b;x_{1},\y\z)Y(v;z_{0},\z)\\
&&-\Res_{x_{1}}y_{0}^{-1}\delta\left(\frac{z_{0}-x_{1}}{-y_{0}}\right)
Y(v;z_{0},\z)Y(Y(a;x_{0},\x)b;x_{1},\y\z).
\end{eqnarray*}
This shows that $a_{m_{0},\m}b\in K$ for $m_{0}\in \Z,\ \m\in
\Z^{r}$. It follows from induction and the span assumption that
$V=K$. This proves that (\ref{eiterate-thm}) holds for {\em all }
$u,v\in V$. Then weak associativity follows. Therefore, $(V,Y)$
carries the structure of an $(r+1)$-toroidal vertex algebra
without vacuum.
\end{proof}

\br{rproof}
{\em Note that in the proof of Theorem \ref{tanalogue},
if we can show that the map $Y$ is injective as for the
one-variable case, the proof can be simplified by using the toroidal
vertex algebra structure transported from $\<\tilde{U}\>$. Unfortunately,
$Y$ may be not injective.}
\er

\section{$(r+1)$-toroidal vertex algebras and modules associated to
toroidal Lie algebras}

In this section,  by using the general construction established in
Section 3 we associate $(r+1)$-toroidal vertex algebras and their
modules to toroidal Lie algebras.

Let $\g$ be a (possibly infinite-dimensional) Lie algebra, equipped
with a (possibly degenerate) symmetric invariant bilinear form
$\<\cdot,\cdot\>$. Form the $(r+1)$-loop Lie algebra
$$L_{r+1}(\g)=\g \otimes \C[t_{0}^{\pm 1},t_{1}^{\pm 1},\dots,t_{r}^{\pm 1}],$$
where $\C[t_{0}^{\pm 1},t_{1}^{\pm 1},\dots,t_{r}^{\pm 1}]$ is the
algebra of Laurent polynomials in mutually commuting variables
$t_{0},t_{1},\dots,t_{r}$. Form a $1$-dimensional central extension
\begin{eqnarray}
\widehat{L_{r+1}(\g)}=L_{r+1}(\g)\oplus \C {\bf k}=\left(\g \otimes
\C[t_{0}^{\pm 1},t_{1}^{\pm 1},\dots,t_{r}^{\pm 1}]\right)\oplus \C
{\bf k},
\end{eqnarray}
where ${\bf k}$ is central and
\begin{eqnarray}\label{eabmn-relation}
[a\otimes t_{0}^{m_{0}}{\bf t}^{\m},b\otimes t_{0}^{n_{0}}{\bf
t}^{\n}]=[a,b]\otimes t_{0}^{m_{0}+n_{0}}{\bf
t}^{\m+\n}+m_{0}\<a,b\>\delta_{m_{0}+n_{0},0}\delta_{\m+\n,0}{\bf k}
\end{eqnarray}
for $a,b\in \g,\; m_{0},n_{0}\in \Z,\; \m,\n\in \Z^{r}$.

Notice that Lie algebra $\widehat{L_{r+1}(\g)}$ can also be
considered as the affine Lie algebra of the $r$-loop Lie algebra
$\g\otimes \C[t_{1}^{\pm 1},\dots,t_{r}^{\pm 1}]$ with respect to
the symmetric invariant bilinear form defined by
$$\<a\otimes \bft^{\m},b\otimes \bft^{\n}\>=\<a,b\>\delta_{\m+\n,0}
\ \ \ \mbox{ for }a,b\in \g,\; \m,\n\in \Z^{r}.$$

For $a\in \g,\;\m\in \Z^{r}$, set
\begin{eqnarray}
a(\m,z)=\sum_{k\in \Z}a(k,\m)z^{-k-1},
\end{eqnarray}
where $a(k,\m)=a\otimes t_{0}^{k}\bft^{\m}$. The Lie bracket
relations (\ref{eabmn-relation}) amount to
\begin{eqnarray}\label{eab-mul}
[a(\m,z_{1}),b(\n,z_{2})] =[a,b](\m+\n,z_{2})
 z_{1}^{-1}\delta\left(\frac{z_{2}}{z_{1}}\right)
+\<a,b\>\delta_{\m+\n,0}{\bf k}\frac{\partial}{\partial
z_{2}}z_{1}^{-1}\delta\left(\frac{z_{2}}{z_{1}}\right)
\end{eqnarray}
for $a,b\in \g,\; \m,\n\in \Z^{r}$. Furthermore, for $a\in \g$, set
\begin{eqnarray}
a(x_{0},\x)=\sum_{n_{0}\in \Z}\sum_{\n\in \Z^{r}}(a\otimes
t_{0}^{n_{0}}{\bf t}^{\n})x_{0}^{-n_{0}-1}\x^{-\n}=\sum_{\n\in
\Z^{r}}a(\n,x_{0})\x^{-\n}.
\end{eqnarray}
Then (\ref{eab-mul}) amounts to
\begin{eqnarray}\label{eambzz}
[a(\m,x_{0}),b(z_{0},\z)]=\z^{\m}\left([a,b](z_{0},\z)
 x_{0}^{-1}\delta\left(\frac{z_{0}}{x_{0}}\right)
+\<a,b\>{\bf k}\frac{\partial}{\partial
z_{0}}x_{0}^{-1}\delta\left(\frac{z_{0}}{x_{0}}\right)\right).
\end{eqnarray}

Define $r+1$ derivations $d_{i}$ $(0\le i\le r)$ on
$\widehat{L_{r+1}(\g)}$ by
\begin{eqnarray}\label{ederivations}
d_{0}=-1\otimes \frac{\partial}{\partial t_{0}},\ \ \ \
d_{i}=-1\otimes t_{i} \frac{\partial}{\partial t_{i}}\ \ \mbox{ for
}1\le i\le r.
\end{eqnarray}
We have
\begin{eqnarray}
[d_{0},a(x_{0},\x)]=\frac{\partial}{\partial x_{0}}a(x_{0},\x),\ \ \
\  [d_{i},a(x_{0},\x)]=x_{i}\frac{\partial}{\partial
x_{i}}a(x_{0},\x)
\end{eqnarray}
for $1\le i\le r,\; a\in \g$. Set
\begin{eqnarray}
\widehat{L_{r+1}(\g)}_{+}&=&\g\otimes \C[t_{0},t_{1}^{\pm
1},\dots,t_{r}^{\pm 1}],\nonumber\\
\widehat{L_{r+1}(\g)}_{-}&=&\g\otimes
t_{0}^{-1}\C[t_{0}^{-1},t_{1}^{\pm 1},\dots,t_{r}^{\pm 1}],
\end{eqnarray}
which are Lie subalgebras of $\widehat{L_{r+1}(\g)}$. Clearly, both
are stable under the actions of derivations
$d_{0},d_{1},\dots,d_{r}$.

We view $\widehat{L_{r+1}(\g)}$ as a $\Z$-graded Lie algebra with
\begin{eqnarray}
\widehat{L_{r+1}(\g)}_{(n)}=
\begin{cases}\g\otimes t_{0}^{-n}\C[t_{1}^{\pm 1},\dots, t_{r}^{\pm 1}]
&\mbox{ if }n\ne 0\\
\g\otimes \C[t_{1}^{\pm 1},\dots, t_{r}^{\pm 1}]+\C {\bf k}&\mbox{
if }n=0
\end{cases}
\end{eqnarray}
for $n\in \Z$. We have
\begin{eqnarray}
d_{0}\cdot \widehat{L_{r+1}(\g)}_{(n)}\subset
\widehat{L_{r+1}(\g)}_{(n+1)} \ \mbox{ and }\ d_{i}\cdot
\widehat{L_{r+1}(\g)}_{(n)}=\widehat{L_{r+1}(\g)}_{(n)}
\end{eqnarray}
for $1\le i\le r$.

\bl{lvacuum-construction} Let $\ell$ be any complex number. Then
there exists an $(\widehat{L_{r+1}(\g)}_{+}+\C {\bf k})$-module
structure on $\g\oplus \C$, which is uniquely determined by ${\bf
k}=\ell $ (a scalar),
\begin{eqnarray}\label{eaction-module}
&&(a\otimes t_{0}^{k}{\bf t}^{\m})\cdot 1=0,\nonumber\\
 &&(a\otimes t_{0}^{k}{\bf t}^{\m})\cdot b
=\begin{cases}[a,b]&\mbox{ if }k=0\\
\<a,b\>\ell &\mbox{ if }k=1\\
0 &\mbox{ if }k\ge 2
\end{cases}
\end{eqnarray}
for $a,b\in \g$, $k\in \Z,\; \m\in \Z^{r}$ with $k\ge 0$.
Furthermore, $\g\oplus \C$ is an $\N$-graded
$(\widehat{L_{r+1}(\g)}_{+}+\C {\bf k})$-module with $\deg \C=0$ and
$\deg \g=1$.\el

\begin{proof} First, we show that $\g\oplus \C$ is a module for Lie algebra
$\g[t_{0}]$ under the action given as a special case with $\m=0$.
This can be proved straightforwardly. Here we give a proof using the
vertex algebra associated to affine Lie algebra $\hat{\g}$. Recall
(see \cite{fz},\; cf. \cite{ll}) that associated to $\hat{\g}$ with
level $\ell$, one has a vertex algebra $V_{\hat{\g}}(\ell,0)$ whose
underlying space is the universal level-$\ell$ vacuum
$\hat{\g}$-module. View $\g\oplus \C$ as a subspace of
$V_{\hat{\g}}(\ell,0)$ by identifying $a+\alpha\in \g+\C$ with
$a(-1){\bf 1}+\alpha {\bf 1}\in V_{\hat{\g}}(\ell,0)$, where ${\bf
1}$ denotes the canonical highest weight vector. For $a,b\in \g,\;
n\ge 0$, we have
$$a(n){\bf 1}=0,\ \ \ \ a(0)b=[a,b],\ \ a(1)b=\ell \<a,b\>{\bf 1},\ \ \
a(m)b=0\ \ \mbox{ for }m\ge 2.$$ It follows that $\g+\C$ is a
$\g[t_{0}]$-submodule of $V_{\hat{\g}}(\ell,0)$ as desired. Second,
equip $\g\oplus \C$ with the evaluation module structure for the
$r$-loop Lie algebra $(\g[t_{0}])\otimes \C[t_{1}^{\pm 1},\dots,t_{r}^{\pm 1}]$ with
${\bf t}=1$ (through the Lie algebra homomorphism sending $a\otimes
t_{0}^{k}\bft^{m}$ to $a\otimes t_{0}^{k}$.) As
$\widehat{L_{r+1}(\g)}_{+}+\C {\bf k}$ is a direct sum of Lie
algebras, by letting ${\bf k}$ act as scalar $\ell$ we obtain an
$(\widehat{L_{r+1}(\g)}_{+}+\C {\bf k})$-module structure on
$\g\oplus \C$, as desired. It is clear that $\g\oplus \C$ is an
$\N$-graded module.
\end{proof}

Let $\ell\in \C$. Denote by $(\g+ \C)_{\ell}$ the
$(\widehat{L_{r+1}(\g)}_{+}+\C {\bf k})$-module obtained in Lemma
\ref{lvacuum-construction}. We form an induced module
\begin{eqnarray}
V_{\widehat{L_{r+1}(\g)}}(\ell,0)
=U\left(\widehat{L_{r+1}(\g)}\right)\otimes_{U(\widehat{L_{r+1}(\g)}_{+}+\C{\bf
k})} (\g+\C)_{\ell},
\end{eqnarray}
which is an $\N$-graded $\widehat{L_{r+1}(\g)}$-module of level
$\ell$. By the P-B-W theorem, we have
$$V_{\widehat{L_{r+1}(\g)}}(\ell,0)=U(\widehat{L_{r+1}(\g)}_{-})\otimes (\g\oplus\C)$$
as a vector space. It follows from the $\N$-grading that
$V_{\widehat{L_{r+1}(\g)}}(\ell,0)$ is a restricted
$\widehat{L_{r+1}(\g)}$-module.

Set $${\bf 1}=1\otimes 1\in V_{\widehat{L_{r+1}(\g)}}(\ell,0).$$
Identify $\g$ with the subspace $1\otimes \g$, through the map
$a\mapsto 1\otimes a$.

We have:

\bt{ttoroidal-algebra} Let $\ell$ be any complex number. Then there
exists an $(r+1)$-toroidal vertex algebra structure on
$V_{\widehat{L_{r+1}(\g)}}(\ell,0)$, which is uniquely determined by
the conditions that ${\bf 1}$ is the vacuum vector and that
$$Y(a;x_{0},\x)=a(x_{0},\x)\ \ \ \mbox{ for }a\in \g.$$
\et

\begin{proof} As
$V_{\widehat{L_{r+1}(\g)}}(\ell,0)=U(\widehat{L_{r+1}(\g)})(\g+\C
{\bf 1})$, the uniqueness follows immediately. We now establish the
existence by applying Theorem \ref{tanalogue}. Let $W$ be any
restricted $\widehat{L_{r+1}(\g)}$-module of level $\ell$. In
particular, we can (and we shall) take $W$ to be
$V_{\widehat{L_{r+1}(\g)}}(\ell,0)$. Set
$$U_{W}=\{ a(x_{0},\x)\;|\; a\in \g\}\cup \{ 1_{W}\}.$$
{}For $a,b\in \g$, from (\ref{eab-mul}) we have
$$(x_{0}-z_{0})^{2}[a(\m,x_{0}),b(\n,z_{0})]=0$$
for $\m,\n\in \Z^{r}$, which implies
$$(x_{0}-z_{0})^{2}[a(x_{0},\x),b(z_{0},\z)]=0.$$
Thus $U_{W}$ is a local subset of $\E(W,r)$. By Theorem
\ref{tclosedness-tva}, $U_{W}$ generates an $(r+1)$-toroidal vertex
algebra $\<U_{W}\>$ inside $\E(W,r)$ and $W$ is a $\<U_{W}\>$-module
with
$$Y_{W}(\beta(x_{0},\x);z_{0},\z)=\beta(z_{0},\z)\ \ \mbox{ for }
\beta(x_{0},\x)\in \<U_{W}\>.$$ For $a,b\in \g,\; \m\in \Z^{r}$,
with (\ref{eambzz}) by Proposition \ref{ptransfer}, we have
\begin{eqnarray*}
&&[Y_{\E}(a(t_{0},\bft); x_{0},\m),Y_{\E}(b(t_{0},\bft);
z_{0},\z)]\\
&=&\z^{\m}\left(Y_{\E}([a,b](t_{0},\bft);z_{0},\z))
 x_{0}^{-1}\delta\left(\frac{z_{0}}{x_{0}}\right)
+\<a,b\>{\bf k}\frac{\partial}{\partial
z_{0}}x_{0}^{-1}\delta\left(\frac{z_{0}}{x_{0}}\right)\right).
\end{eqnarray*}
Thus, $\<U_{W}\>$ becomes an $\widehat{L_{r+1}(\g)}$-module of level
$\ell$ with $a(z_{0},\z)$ acting as $Y_{\E}(a(x_{0},\x);z_{0},\z)$
for $a\in \g$.  Let $a,b\in \g,\; m_{0}\in \Z,\; \m\in \Z^{r}$ with
$m_{0}\ge 0$. With relation (\ref{eambzz}), in view of Lemma
\ref{ltransfer} we have
$$a(x_{0},\x)_{m_{0},\m}1_{W}=0,$$
and
\begin{eqnarray*}
a(x_{0},\x)_{m_{0},\m}b(x_{0},\x)=\begin{cases}[a,b](x_{0},\x)&
\mbox{ for }m_{0}=0\\
\ell \<a,b\> 1_{W} &\mbox{ for }m_{0}=1\\
0& \mbox{ for }m_{0}\ge 2.
\end{cases}
\end{eqnarray*}
This implies that $U_{W}+\C 1_{W}$ is an
$(\widehat{L_{r+1}(\g)}_{+}+\C {\bf k})$-submodule of $\<U_{W}\>$
and that the obvious map $\psi_{W}$ {}from $\g+ \C$ to $U_{W}+\C
1_{W}$ is an $(\widehat{L_{r+1}(\g)}_{+}+\C {\bf k})$-module
homomorphism. It then follows from the construction of
$V_{\widehat{L_{r+1}(\g)}}(\ell,0)$ that there exists an
$\widehat{L_{r+1}(\g)}$-module homomorphism $\tilde{\psi}_{W}$ from
$V_{\widehat{L_{r+1}(\g)}}(\ell,0)$ to $\<U_{W}\>$, extending
$\psi_{W}$.

Now, take $W=V_{\widehat{L_{r+1}(\g)}}(\ell,0)$. Denote
$\tilde{\psi_{W}}$ by $\tilde{\psi}_{x_{0},\x}$ to indicate the
dependence on those variables and also denote
$1_{V_{\widehat{L_{r+1}(\g)}}(\ell,0)}$ simply by $1$. For $v\in
V_{\widehat{L_{r+1}(\g)}}(\ell,0)$, set
$$Y(v;x_{0},\x)=\tilde{\psi}_{x_{0},\x}(v). $$ For $a\in \g$, we have
$$Y(a;x_{0},\x)=\tilde{\psi}_{x_{0},\x}(a)=a(x_{0},\x).$$
Then $\{ Y(a;x_{0},\x)\;|\; a\in \g\}\cup \{1\}$ is local. For $a\in
\g,\ v\in V_{\widehat{L_{r+1}(\g)}}(\ell,0)$, we have
\begin{eqnarray*}
&&Y(Y(a;z_{0},\z)v;y_{0},\y)=\tilde{\psi}_{y_{0},\y}(a(z_{0},\z)v)\\
&=&Y_{\E}(a(y_{0},\y);z_{0},\z)\tilde{\psi}_{y_{0},\y}(v)\\
&=&\Res_{x_{0}}z_{0}^{-1}\delta\left(\frac{x_{0}-y_{0}}{z_{0}}\right)
 a(x_{0},\z\y)\tilde{\psi}_{y_{0},\y}(v)-z_{0}^{-1}\delta\left(\frac{y_{0}-x_{0}}{-z_{0}}\right)
 \tilde{\psi}_{y_{0},\y}(v)a(x_{0},\z\y)\\
&=&\Res_{x_{0}}z_{0}^{-1}\delta\left(\frac{x_{0}-y_{0}}{z_{0}}\right)
Y(a;x_{0},\z\y)Y(v;y_{0},\y)\nonumber\\
& &\hspace{.5cm} -\Res_{x_{0}}
z_{0}^{-1}\delta\left(\frac{y_{0}-x_{0}}{-z_{0}}\right)
Y(v;y_{0},\y)Y(a;x_{0},\z\y).
\end{eqnarray*}
Now it follows immediately from Theorem \ref{tanalogue} that
$(V_{\widehat{L_{r+1}(\g)}}(\ell,0),Y)$ carries the structure of an
$(r+1)$-toroidal vertex algebra as
desired.
\end{proof}

Furthermore, we have:

\bt{tmodule-appl} Let $W$ be any restricted
$\widehat{L_{r+1}(\g)}$-module of level $\ell$. Then there exists a
$V_{\widehat{L_{r+1}(\g)}}(\ell,0)$-module structure on $W$, which
is uniquely determined by
$$Y_{W}(a;x_{0},\x)=a(x_{0},\x)\ \ \ \mbox{ for }a\in \g.$$
On the other hand, suppose that $(W,Y_{W})$ is a
$V_{\widehat{L_{r+1}(\g)}}(\ell,0)$-module. Then $W$ becomes a
restricted $\widehat{L_{r+1}(\g)} $-module of level $\ell$ with
$$a(x_{0},\x)=Y_{W}(a;x_{0},\x)\ \ \ \mbox{ for }a\in \g.$$
 \et

\begin{proof} Recall the first part of the proof of Theorem
\ref{ttoroidal-algebra}: We have an $(r+1)$-toroidal vertex algebra
$\<U_{W}\>$ with $W$ as a canonical module and we showed that
$\<U_{W}\>$ is an $\widehat{L_{r+1}(\g)}$-module with $a(z_{0},\z)$
acting as $Y_{\E}(a(x_{0},\x);z_{0},\z)$ for $a\in \g$ and that
there exists an $\widehat{L_{r+1}(\g)}$-module homomorphism
$\tilde{\psi}_{W}$ from $V_{\widehat{L_{r+1}(\g)}}(\ell,0)$ to
$\<U_{W}\>$, sending ${\bf 1}$ to $1_{W}$ and $a\in \g$ to
$a(x_{0},\x)$.  For $a\in \g,\; v\in
V_{\widehat{L_{r+1}(\g)}}(\ell,0)$, we have
$$\tilde{\psi}_{W}(Y(a;y_{0},\y)v)=\tilde{\psi}_{W}(a(y_{0},\y)v)
=Y_{\E}(a(x_{0},\x);y_{0},\y)\tilde{\psi}_{W}(v)
=Y_{\E}(\tilde{\psi}_{W}(a);y_{0},\y)\tilde{\psi}_{W}(v).$$ Since
$\g$ generates $V_{\widehat{L_{r+1}(\g)}}(\ell,0)$, it follows that
$\tilde{\psi}_{W}$ is a homomorphism of $(r+1)$-toroidal vertex
algebras. With $W$ as a canonical $\<U_{W}\>$-module, consequently,
$W$ becomes a $V_{\widehat{L_{r+1}(\g)}}(\ell,0)$-module with
$$Y_{W}(v;z_{0},\z)=Y_{W}(\tilde{\psi}_{W}(v);z_{0},\z)=\tilde{\psi}_{W}(v)(z_{0},\z)$$
for $v\in V_{\widehat{L_{r+1}(\g)}}(\ell,0)$. In particular, we have
$$Y_{W}(a;z_{0},\z)=Y_{W}(\tilde{\psi}_{W}(a);z_{0},\z)
=Y_{W}(a(x_{0},\x);z_{0},\z)=a(z_{0},\z)$$
for $a\in \g$. This proves the first assertion.

On the other hand, let $(W,Y_{W})$ be a
$V_{\widehat{L_{r+1}(\g)}}(\ell,0)$-module. {}From
(\ref{eaction-module}) we have
\begin{eqnarray}
&&a_{k,\m}\cdot 1=0,\nonumber\\
 &&a_{k,\m}\cdot b
=\begin{cases}[a,b]&\mbox{ if }k=0\\
\<a,b\>\ell{\bf 1} &\mbox{ if }k=1\\
0 &\mbox{ if }k\ge 2
\end{cases}
\end{eqnarray}
for $a,b\in \g$, $k\ge 0,\; \m\in \Z^{r}$. Combining this with the
commutator formula (\ref{eumv}), we get
\begin{eqnarray*}
&&\Res_{\x}\x^{\m-1}[Y_{W}(a;x_{0},\x),Y_{W}(b;z_{0},\z)]\\
&=&\z^{\m}\left(Y_{W}([a,b];z_{0},\z)x_{0}^{-1}\delta\left(\frac{z_{0}}{x_{0}}\right)+\ell \<a,b\>
Y_{W}({\bf 1};z_{0},\z)\frac{\partial}{\partial
z_{0}}x_{0}^{-1}\delta\left(\frac{z_{0}}{x_{0}}\right)\right).
\end{eqnarray*}
It follows that $W$ is a restricted $\widehat{L_{r+1}(\g)}$-module
of level $\ell$ with $a(x_{0},\x)=Y_{W}(a;x_{0},\x)$ for $a\in \g$.
\end{proof}

\br{rkernel} {\em Recall from Proposition \ref{pva-subalgebra} that
for any $(r+1)$-toroidal vertex algebra $V$, the left ideal $V^{0}$ generated by ${\bf 1}$
 is an ordinary vertex algebra, where
\begin{eqnarray*}
V^{0}&=&\{ u_{m_{0},\m}{\bf 1}\;|\; u\in V,\; m_{0}\in \Z,\; \m\in
\Z^{r}\},
\end{eqnarray*} and
$$Y^{0}(v,x_{0})=Y(v;x_{0},\x)|_{\x=1} \ \ \ \mbox{ for }v\in V^{0}.$$
Assume $V=V_{\widehat{L_{r+1}(\g)}}(\ell,0)$.  As $\g+\C {\bf 1}$ generates $V_{\widehat{L_{r+1}(\g)}}(\ell,0)$,
it follows that $$V^{0}=U(\widehat{L_{r+1}(\g)}){\bf 1}=V_{\widehat{L_{r}(\g)}}(\ell,0)$$
as a vector space, where $V_{\widehat{L_{r}(\g)}}(\ell,0)$ is the vertex algebra associated to the
affine Lie algebra $\widehat{L_{r}(\g)}$ of level $\ell$.
For $a\in \g,\; \m\in \Z^{r}$, using (\ref{e-mid-step}) we have
$$Y^{0}((a\otimes t_{0}^{-1}\otimes {\bf t}^{\m}){\bf 1},x_{0})
=Y((a\otimes t_{0}^{-1}\otimes {\bf t}^{\m}){\bf 1};x_{0},\x)|_{\x=1}
=Y(a;x_{0},\m)=a(x_{0},\m).$$
It follows (from induction) that $V^{0}$ coincides with the vertex algebra $V_{\widehat{L_{r}(\g)}}(\ell,0)$.} \er

In what follows, by modifying the previous construction we shall
construct an extended $(r+1)$-toroidal vertex algebra. First, we
extend the toroidal Lie algebra $\widehat{L_{r+1}(\g)}$ by
derivations $d_{0},d_{1},\dots,d_{r}$ as defined in
(\ref{ederivations}).
 Set
\begin{eqnarray}\label{eDspace}
D=\C d_{0}+\C d_{1}+\cdots +\C d_{r},
\end{eqnarray}
an abelian Lie algebra acting on $\widehat{L_{r+1}(\g)}$ by
derivations. We have a Lie algebra
\begin{eqnarray}
\widetilde{L_{r+1}(\g)}=\widehat{L_{r+1}(\g)}\oplus D,
\end{eqnarray}
which is the semi-direct product Lie algebra. Recall  the
$(\widehat{L_{r+1}(\g)}_{+}+\C{\bf k})$-module $(\g+\C)_{\ell}$ with
$\ell\in \C$ (see Lemma \ref{lvacuum-construction}). Form an induced
$(\widehat{L_{r+1}(\g)}_{+}+\C{\bf k}+D)$-module
\begin{eqnarray*}
\tilde{K}_{\ell}=U\left(\widehat{L_{r+1}(\g)}_{+}+\C{\bf
k}+D\right)\otimes_{U(\widehat{L_{r+1}(\g)}_{+}+\C{\bf k})}
(\g+\C)_{\ell}.
\end{eqnarray*}
It can be readily seen that $U(D)D\cdot 1$ is a submodule of
$\tilde{K}_{\ell}$. Define $K_{\ell}$ to be the quotient module of
$\tilde{K}_{\ell}$ modulo relation $D\cdot 1=0$. As a $D$-module,
\begin{eqnarray}\label{eKspace}
K_{\ell}=\left(S(D)\otimes \g\right)\oplus \C,
\end{eqnarray}
the direct sum of the free $D$-module on $\g$ and the trivial module
$\C$.

Then, we form an induced $\widetilde{L_{r+1}(\g)} $-module
\begin{eqnarray}
V_{\widetilde{L_{r+1}(\g)}}(\ell,0)
=U\left(\widetilde{L_{r+1}(\g)}\right)\otimes_{U(\widehat{L_{r+1}(\g)}_{+}+D+\C{\bf
k})} K_{\ell},
\end{eqnarray}
which is an $\N$-graded $\widetilde{L_{r+1}(\g)}$-module of level
$\ell$. We have:

\bt{textended-tva} On $V_{\widetilde{L_{r+1}(\g)}}(\ell,0)$, there
exists a structure of an extended $(r+1)$-toroidal vertex algebra,
which is uniquely determined by the condition
$$\D_{i}=d_{i}\ \mbox{ for }0\le i\le r\ \mbox{ and }\ \
Y(a;x_{0},\x)=a(x_{0},\x)\ \mbox{ for }a\in \g.$$ \et

\begin{proof} First we shall use Theorem \ref{tanalogue} to obtain an $(r+1)$-toroidal vertex
algebra structure. Recall
$$K_{\ell}=\C[d_{0},d_{1},\dots,d_{r}]\g+ \C \subset V_{\widetilde{L_{r+1}(\g)}}(\ell,0).$$
Let $W$ be any restricted $\widetilde{L_{r+1}(\g)}$-module of level
$\ell$. Set
$$U_{W}=\{ d_{i}^{k}a(x_{0},\x)\ |\ a\in \g,\ 0\le i\le r,\ k\ge 0\}.$$
Then $U_{W}$ is a local subset of $\E(W,r)$. By Theorem
\ref{tclosedness-tva}, $U_{W}$ generates an extended
$(r+1)$-toroidal vertex algebra
$\C[D_{0},D_{1},\dots,D_{r}]\<U_{W}\>$ in $\E(W,r)$, denoted locally
by $\overline{\<U_{W}\>}$. It follows from the same argument as in
the proof of Theorem \ref{ttoroidal-algebra} that
$\overline{\<U_{W}\>}$ is an $\widehat{L_{r+1}(\g)}$-module of level
$\ell$ with
$$a(z_{0},\z)=Y_{\E}(a(x_{0},\x);z_{0},\z)\ \ \mbox{ for }a\in \g.$$
As $\overline{\<U_{W}\>}$ is an extended $(r+1)$-toroidal vertex
algebra with $\D_{i}=D_{i}$ $(i=0,1,\dots,r)$, we have
\begin{eqnarray*}
&&[D_{0},Y_{\E}(a(x_{0},\x);z_{0},\z)]=\frac{\partial}{\partial
z_{0}}Y_{\E}(a(x_{0},\x);z_{0},\z),\\
&&[D_{i},Y_{\E}(a(x_{0},\x);z_{0},\z)]=z_{i}\frac{\partial}{\partial
z_{i}}Y_{\E}(a(x_{0},\x);z_{0},\z)
\end{eqnarray*}
for $a\in \g,\; 1\le i\le r$. Thus $\overline{\<U_{W}\>}$ is an
$\widetilde{L_{r+1}(\g)}$-module of level $\ell$ with $d_{i}=D_{i}$
for $0\le i\le r$. From the proof of Theorem \ref{ttoroidal-algebra}
we see that the linear map $\psi_{W}: \g+\C\rightarrow
\overline{\<U_{W}\>}$, defined by
$$\psi_{W} (a+\lambda)=a(x_{0},\x)+\lambda 1_{W}
\ \ \mbox{ for }a\in \g,\; \lambda\in \C,$$ is  an
$(\widehat{L_{r+1}(\g)}_{+}+\C {\bf k})$-module homomorphism. We
also have
$$d_{i}\cdot 1_{W}=D_{i}(1_{W})=0\ \ \mbox{ for }0\le i\le r.$$
 It follows from the construction of
$V_{\widetilde{L_{r+1}(\g)}}(\ell,0)$ that there exists an
$\widetilde{L_{r+1}(\g)}$-module homomorphism $\widetilde{\psi_{W}}$
from $V_{\widetilde{L_{r+1}(\g)}}(\ell,0)$ to
$\overline{\<U_{W}\>}$, extending $\psi_{W}$.

Take $W=V_{\widetilde{L_{r+1}(\g)}}(\ell,0)$ and let
$\widetilde{\psi}(x_{0},\x)$ denote the corresponding map
$\widetilde{\psi_{W}}$. Define
$$Y(v;x_{0},\x)=\widetilde{\psi}(x_{0},\x)(v)\ \ \mbox{ for }
v\in V_{\widetilde{L_{r+1}(\g)}}(\ell,0).$$ We have
\begin{eqnarray*}
&&Y(d_{0}^{k}a;x_{0},\x)=\widetilde{\psi}(x_{0},\x)(d_{0}^{k}a)
=D_{0}^{k}\left(\widetilde{\psi}(x_{0},\x)(a)\right)
=\left(\frac{\partial}{\partial x_{0}}\right)^{k}Y(a;x_{0},\x),\\
&&Y(d_{i}^{k}a;x_{0},\x)=\widetilde{\psi}(x_{0},\x)(d_{i}^{k}a)
=D_{i}^{k}\left(\widetilde{\psi}(x_{0},\x)(a)\right)
=\left(x_{i}\frac{\partial}{\partial x_{i}}\right)^{k}Y(a;x_{0},\x)
\end{eqnarray*}
for $a\in \g,\ k\ge 0,\ 1\le i\le r$. As in the proof of Theorem
\ref{ttoroidal-algebra}, it follows from Theorem \ref{tanalogue}
with $U=K_{\ell}$ that $(V_{\widetilde{L_{r+1}(\g)}}(\ell,0),Y)$
carries the structure of an $(r+1)$-toroidal vertex algebra. We also
have
\begin{eqnarray*}
&&[d_{0},Y(a;x_{0},\x)]=[d_{0},a(x_{0},\x)]=\frac{\partial}{\partial
x_{0}}a(x_{0},\x)=\frac{\partial}{\partial
x_{0}}Y(a;x_{0},\x),\\
&&[d_{i},Y(a;x_{0},\x)]=[d_{i},a(x_{0},\x)]=x_{i}\frac{\partial}{\partial
x_{i}}a(x_{0},\x)=x_{i}\frac{\partial}{\partial x_{i}}Y(a;x_{0},\x).
\end{eqnarray*}
Then it follows that
\begin{eqnarray*}
&&[d_{0},Y(d_{j}^{k}a;x_{0},\x)]=\frac{\partial}{\partial
x_{0}}Y(d_{j}^{k}a;x_{0},\x),\\
&&[d_{i},Y(d_{j}^{k}a;x_{0},\x)]=x_{i}\frac{\partial}{\partial
x_{i}}Y(d_{j}^{k}a;x_{0},\x).
\end{eqnarray*}
 As $K_{\ell}$ generates $V_{\widetilde{L_{r+1}(\g)}}(\ell,0)$,
by Lemma \ref{lextended-generating},
$V_{\widetilde{L_{r+1}(\g)}}(\ell,0)$ is an extended
$(r+1)$-toroidal vertex algebra with $\D_{i}=d_{i}$
$(i=0,1,\dots,r)$.
\end{proof}

We also have:

\bt{textended-module} For any restricted
$\widetilde{L_{r+1}(\g)}$-module $W$ of level $\ell$, there exists a
structure of a $V_{\widetilde{L_{r+1}(\g)}}(\ell,0)$-module $Y_{W}$
on $W$, uniquely determined by
$$Y_{W}(a;x_{0},\x)=a(x_{0},\x)\ \ \mbox{ for }a\in \g.$$
On the other hand, for any
$V_{\widetilde{L_{r+1}(\g)}}(\ell,0)$-module $(W,Y_{W})$, $W$ is a
restricted $\widetilde{L_{r+1}(\g)}$-module of level $\ell$ with
$$a(x_{0},\x)=Y_{W}(a;x_{0},\x)\ \ \mbox{ for }a\in \g\ \ \mbox{ and }\
d_{i}=D_{i}\ \ \mbox{ for }0\le i\le r.$$
 \et

\br{rmodule-hal-D} {\em Let $W$ be a restricted
$\widehat{L_{r+1}(\g)}$-module of level $\ell$. Then using Theorem
\ref{tclosedness-tva} one can show that $W$ is a module for
$V_{\widetilde{L_{r+1}(\g)}}(\ell,0)$ viewed as an $(r+1)$-toroidal
vertex algebra, satisfying
\begin{eqnarray*}
&&Y_{W}(\D_{0}(v);x_{0},\x)=\frac{\partial}{\partial
x_{0}}Y_{W}(v;x_{0},\x),\\
&&Y_{W}(\D_{i}(v);x_{0},\x)=x_{i}\frac{\partial}{\partial
x_{i}}Y_{W}(v;x_{0},\x)
\end{eqnarray*}
for $v\in V,\; 1\le i\le r$.} \er

\end{document}